\documentclass[pdflatex,sn-mathphys-num]{sn-jnl}

\usepackage{graphicx}%
\usepackage{multirow}%
\usepackage{amsmath,amssymb,amsfonts}%
\usepackage{amsthm}%
\usepackage{mathrsfs}%
\usepackage[title]{appendix}%
\usepackage{xcolor}%
\usepackage{textcomp}%
\usepackage{manyfoot}%
\usepackage{booktabs}%
\usepackage{algorithm}%
\usepackage{algorithmicx}%
\usepackage{algpseudocode}%
\usepackage{listings}%

\usepackage{psfrag} 
\usepackage{subfig} 

\theoremstyle{thmstyleone}%
\theoremstyle{thmstyletwo}%

\theoremstyle{thmstylethree}%

\begin{document}

\title[PO-CSA:title]{Perpertual Coupled Simulated Annealing for Continuous Optimization}


\author*[1,3]{\fnm{Kayo} \sur{Goncalves-e-Silva}}\email{kayo@imd.ufrn.br}
\equalcont{These authors contributed equally to this work.}

\author[2,3]{\fnm{Samuel} \sur{Xavier-de-Souza}}\email{samuel@dca.ufrn.br}
\equalcont{These authors contributed equally to this work.}


\affil*[1]{\orgdiv{Digital Metropolis Institute}}


\affil[2]{\orgdiv{Department of Computation and Automation}}

\affil[3]{ \orgname{Federal University of Rio Grande do Norte}, \orgaddress{\street{Campus Universitário Lagoa Nova}, \city{Natal}, \state{Rio Grande do Norte}, \postcode{59078-900}, \country{Brazil}} }



\abstract{Global optimization heuristics are popular to optimize hard non-convex problems. Despite their irrefutably large cost-to-solution, in the lack of other working greedy or convex approaches, global optimization algorithms remain the no-brainer choice. Nevertheless, successful use often requires tedious adjustments to initial parameters to avoid premature convergence. The Coupled Simulated Annealing approach proposed a method based on the coupling of multiple optimizers to escape premature convergence, having achieved success in optimizing hyperparameters of many applications of machine learning; however, a careful choice of the generation temperature is still required. In this paper we propose the Perpetual Orbit technique as a solution
to control the generation temperature and avoid search stagnation. In
principle, this technique can also be applied to other ensemble- and population-based algorithms that have a dispersion variable. The results of our experiments show superior performance when using the proposed technique because it makes the Couple Simulated Annealing totally parameter-free and capable of reaching equal or better solutions in more than 85\% cases across all functions and competitor methods analyzed.}

\keywords{Coupled Simulated Annealing, Perpetual Orbit, Global Optimization, Parameter-free Optimization}



\maketitle

\section{Introduction}\label{sec1}

The research interest in Metaheuristics for global optimization have been increasing along the past decades, becoming a very fertile field for new strategies~\citep{pan2012new,das2014spatially,mahdavi2015metaheuristics,mirjalili2015ant,zhang2016memetic,zhang2017novel,hussain2019metaheuristic,halim2021performance,ribeiro2024,juan2023review}. Among other positive characteristics, metaheuristics are usually simple to implement, do not require gradient information, can scape from local optima and can be applied to a wide range of real-world optimization~\citep{ding2015transaction,mirjalili2016whale,mistry2017microga,osaba2021tutorial,marti2025fifty,juan2023review}. 

Several established  algorithms for global optimization require a stochastic distribution to generate new solutions. It is usually necessary to control a variable responsible for the dispersion of the distribution. In Simulated Annealing (SA)~\citep{Kirkpatrick83optimizationby}, for example, the generation and acceptance temperatures are dispersion variables; in Genetic Algorithms, it is a factor in the Blend Crossover~\citep{eshelman1993real} and in the Linear Crossover~\citep{wright1991genetic}; and, in the Swarm Particle~\citep{kenny1995particle}, it is the cognitive and social components. Such variables are hard to tune for several reasons, including: the lack of statistical information; the poorly defined operating range; the need to initialize and schedule them; the dependency on the optimization problem; and the demand for a very expensive tuning process due to its empirical nature.

The Coupled Simulated Annealing (CSA) algorithm~\citep{samuel2010} is a method for global optimization of continuous variables based on SA. It is composed by several distributed SA processes with coupled acceptance probabilities. It is mainly characterized by its acceptance probability function and a coupling term, which is a function of the current energies of all SA processes. CSA is more robust than the uncoupled case with respect to initialization parameter because the acceptance temperature does not require tuning. However, the choice of the generation temperature still affects the optimization process, and defining an initial value and a monotonic annealing scheduling is necessary.

In this paper, we propose a novel technique called Perpetual Orbit (PO) that can control variables responsible for the dispersion of the distribution without prior informations. The main purpose is to maintain or improve the quality of the final solution with an automatic parameter control that, in addition to  avoid the convergence to local optima,  facilitates the application of metaheuristics by reducing the dependence on metaheuristic's success rate on the tuning of its initial parameters. In order to prove the effectiveness of the PO technique, we applied it to the generation temperature of the CSA algorithm. The resulting algorithm, called PO-CSA, is a parameter-free algorithm that avoids stagnation of the CSA's generation temperature. We performed experiments to compare the proposed PO-CSA with classic metaheuristics, such as Cuckoo Search~\citep{yang2009cuckoo}, Genetic Algorithm~\citep{holland1992adaptation}, Differential Evolution~\citep{storn1997differential} and Particle Swarm Optimization~\citep{kenny1995particle}. We also compared the PO-CSA to the original CSA tuned with an exhaustive search of the initial generation temperature. The results show that the PO-CSA is robust with respect to the initial parameters. It obtains better solutions in most of the cases and outperforms the CSA even with the initial parameters adjusted exhaustively.


\section{Coupled Simulated Annealing}
\label{sec:csa}


Coupled Simulated Annealing (CSA)~\cite{samuel2010} is a stochastic global optimization method capable of reducing the sensitivity of the initialization parameters while guiding the optimization of continuous variables to quasi-optimal runs. It is based on the Simulated Annealing (SA)~\cite{Kirkpatrick83optimizationby} and the Coupled Local Minimizers (CLM)~\cite{suykens2001intelligence} algorithms.

The CSA algorithm consists of an ensemble of SA optimizers whose behavior, individually, is similar to the execution of a SA process, i.e. the algorithmic steps that involve generation and acceptance of a single current solution are separately performed for each SA optimizer. In the generation step,  the generation temperature $T^{\rm gen}_k$ is responsible for the degree of similarity between the current solution and the probing solution at iteration $k$, i.e. the larger the value of the generation temperature, the larger the average difference between current and probing solutions. The acceptance temperature $T^{\rm ac}_k$, used in the acceptance step at iteration $k$, composes the calculation of the chance of accepting worse probing solutions, i.e., the larger the value of the acceptance temperature, the better the chances of worse probing solutions be accepted.


Unlike SA, where the acceptance probability is a scalar function of the current and probing solutions alone, the CSA acceptance probability function is a scalar function according to
\begin{equation}
	0 \leq A_\Theta (\gamma, x_i \to y_i ) \leq 1,
\end{equation}
for every $x_i \in \Theta$, $y_i \in \Omega$, and $\Omega \in \mathbb{R}^m$, with $\Theta$ and $\Omega$ being the set of current and the set of all solutions, respectively, and $i=1,\ldots,m$, with $m$ as the number of elements of $\Theta$. Each solutions $x_i$ has an associated energy $E(x_i)$, equivalent to the cost of the solution $x_i$. Finally, the term $\gamma$ is responsible for the coupling between the SA optimizers. 

The acceptance probability function used in this work avoids numerical instabilities due to the energy of the states being positive or negative. It is defined by:
\begin{equation}
\label{eq:acceptance_prob_function}
\displaystyle
A_\Theta\left(\gamma, x_i \to y_i \right) = \frac{\displaystyle \exp\left( \frac{E(x_i)-max\left(E\left(x_i\right)\right)_{x_i\in \Theta}  }{\displaystyle  T^{\rm ac}_k } \right) }{\gamma},
\end{equation}
where $T^{\rm ac}_k$ is the acceptance temperature at iteration $k$ and $\gamma$ is the coupling term, defined by:
\begin{equation}
\label{eq:coupling_term}
\gamma = \sum_{\forall x \in \Theta}^{}\exp\left( \frac{E(x)-max\left(E\left(x\right)\right)_{x_i\in \Theta}  }{  T^{\rm ac}_k  } \right).
\end{equation}

As the variance of $A_\Theta$ is limited by
\begin{equation}
 0 \le \sigma^2 \le \frac{m-1}{m^2},
\end{equation}
a simple rule can be applied to control this variance: 
\begin{equation}
    \label{eq:cond_var}
    \begin{split}
    \mbox{if} \quad \sigma^2 <   \sigma^2_D, \quad \quad T^{\rm ac}_{k+1} = T^{\rm ac}_{k}\left(1-\alpha \right),\\
    \mbox{if} \quad \sigma^2 \ge \sigma^2_D, \quad \quad T^{\rm ac}_{k+1} = T^{\rm ac}_{k}\left(1+\alpha \right),
    \end{split}
\end{equation}
where $\sigma^2_D$ is the desired variance and $\alpha$ the deprecation rate of the temperature, usually values between $(0.0,0.1]$~\cite{samuel2010}. The acceptance temperature is decreased by a factor of $(1-\alpha)$ if the acceptance variance is below the desired variance, and increased by a factor of $(1+\alpha)$ otherwise. It has been shown that the optimization performance of CSA improves to a quasi-optimal run when the variance of $A_\Theta$ is around $99\%$ of its maximum value~\cite{samuel2010}. This way, the variance control substitutes the monotonic schedule of the acceptance temperature.
This is very important because it avoids the initial setup of the acceptance temperature in classical SA that often demands an exhaustive search. In contrast, two other initial parameters are introduced, $\alpha$ and $\sigma^2_D$, but these parameters have a well-defined operating range and are much less dependent on the optimization problem at hand~\cite{samuel2010}.

The probing solutions are generated through the Cauchy Distribution according to:
\begin{equation}
	y_i=x_i+\epsilon*T_k^{\rm gen},
	\label{eq:cauchy_d}
\end{equation}
where each element of $\epsilon \in \mathbb{R}^m$ is a random independent variable sampled from the Cauchy distribution.

The pseudo-code of the CSA is described in Algorithm~\ref{alg:CSA}.

\begin{algorithm}
\caption{Coupled Simulated Annealing Algorithm}\label{alg:CSA}
\begin{algorithmic}[1]
\State \textbf{Step 1) Initialization:}
        \State \quad - Assign $m$ random initial solutions to $\Theta$; 
        \State \quad - Assign $m$ random initial solutions to $\Theta$;
	\State \quad - Assess the energies $E(x_i)$, $\forall x_i \in \Theta$;
	\State \quad - Set the iteration number $k=0$;
	\State \quad - Set the initial temperatures $T^{\rm gen}_k$ and $T^{\rm ac}_k$;
	\State \quad - Assign random value to $\alpha$; 
	\State \quad - Evaluate the coupling term $\gamma$; 
	\State \quad - Set $\sigma_D^2 = 0.99\frac{(m-1)}{m^2}$. 

\item[]

\State \textbf{Step 2) Generation:}
	\State \quad - Generate a probing solution $y_i=x_i+\epsilon*T_k^{\rm gen}$, $\forall i=1,\cdots,m$,  according to (\ref{eq:cauchy_d}); 
	\State \quad - Assess the energies $E(y_i), \forall i = 1, ~ \ldots ~ , m$.

\item[]

\State \textbf{Step 3) Acceptance:}\;
	\State \quad - $\forall i=1,\cdots,m$, do $x_i \gets y_i$ if: \;
	\State \quad \quad \quad A) $E(y_i) \le E(x_i)$; or \;
	\State \quad \quad \quad B) $A_\Theta > r$, where $r$ is a random variable sampled from a uniform distribution $[0,1]$. \;

\item[]

\State \textbf{Step 4) Update:}\;
	\State \quad - Calculate $\sigma^2$; \;
	\State \quad - Adjust $T^{\rm ac}_k$ according to the following rules: \;
	\State \quad \quad \quad - If $\sigma^2 < \sigma_D^2$ then $T^{\rm ac}_{k+1} = T^{\rm ac}_k(1-\alpha)$; \;
	\State \quad \quad \quad - If $\sigma^2 > \sigma_D^2$ then $T^{\rm ac}_{k+1} = T^{\rm ac}_k(1+\alpha)$. \;
	\State \quad - Decrease the generation temperature $T^{\rm gen}_k$ according to a chosen monotonic schedule; \;
    \State \quad - Evaluate the coupling term $\gamma$; \;
    \State \quad - Increment $k$ by $1$. \;

\State \textbf{Step 5) Stopping Criterion:}\;
	\State \quad - Stop if the stopping criterion is met; \;
	\State \quad - Otherwise, go to Step 2. \;

\end{algorithmic}
\end{algorithm}


\section{Perpetual Orbit Technique}
\label{sec:pot}

The Perpetual Orbit (PO) technique, proposed in this work, is characterized by an ensemble of solutions $s_i \in \Theta$ that cooperate in the effort to find the best value for the dispersion variable during each moment of the execution of a chosen algorithm. 

The main focus of this technique is to maintain or improve the quality of the final solution with an automatic parameter control that drives the optimization resiliently to avoid convergence to local optima and reduces the need for initial algorithm settings based on the problem to be applied. It can be potentially applied to other ensemble or population-based approaches. For that, the technique requires that each solution $s_i$ has its own dispersion variable $V^{s_i}_k$ at iteration $k$. 
Initially, each solution $s_i$ and dispersion variable $V^{s_i}_k$ are randomly initialized and the best initial solution is determined by assessing the costs of all solutions in the ensemble. The value of each $V^{s_i}_k$ changes between its upper and lower bounds, $U^{s_i}_k$ and $L^{s_i}_k$, respectively, derived from the best solution's dispersion variable as follows.
\begin{equation}
	\begin{split}
		U^{s_i}_k = \beta *V_{k}^{s_{\rm best}}\\
		L^{s_i}_k = \frac{1}{\beta} * V_{k}^{s_{\rm best}}\\
	\end{split}
	\label{eq:U_L_initial_bounds}
\end{equation}
where $s_{\rm best}$ is the best solution in $\Theta$, and $\beta$ is the boundary multiplier, a valued in the range $[1,100]$. The parameter $\beta$ is not mandatory and can be removed from the code by simply doing $\beta$ equals 1.

During the optimization, the value of all dispersion variables moves around the value of the best solution's dispersion variable, except for itself, which remains unchanged. At each iteration, each dispersion value moves towards either the upper or lower bound. When reaching any of these bounds, it changes its direction and moves towards the other bound. If the upper bound is reached, the bound value is multiplied by a factor of $(1+\mu)$, for $0 < \mu < 1$. If it is the lower bound that is reached, it is multiplied by $(1-\mu)$.

The movements of increment and decrement on the value of the dispersion variable, as well as the changes on the bounds, occur continuously until the algorithm finds a new best solution. Then, the new best solution's dispersion variable stops changing its values, and the old best solution's dispersion value begins to change around the new value. The upper and lower bound are then reinitialized according to~(\ref{eq:U_L_initial_bounds}).

\begin{figure}[!ht]
	\centering
	\psfrag{a}[][][1.0]{$V^{s_i}_k$} 
	\psfrag{b}[][][1.0]{$~~s_i$} 
	\psfrag{c}[][][1.0]{$~~s_1$}
	\psfrag{d}[][][1.0]{$~~s_2$}
	\psfrag{e}[][][1.0]{$~~s_3$}
	\psfrag{f}[][][1.0]{$~~s_4$}
	\psfrag{g}[][][1.0]{$~~s_5$}
	\psfrag{h}[][][1.0]{$L^{s_i}_k$}
	\psfrag{i}[][][1.0]{$D^{s_2}_k=-1$}
	\psfrag{j}[][][1.0]{$D^{s_1}_k=-1$}
	\psfrag{k}[][][1.0]{$D^{s_3}_k=0$}
	\psfrag{l}[][][1.0]{$D^{s_4}_k=1$}
	\psfrag{m}[][][1.0]{$D^{s_5}_k=1$}
	\psfrag{n}[][][1.0]{$best~s_i$}
	\psfrag{o}[][][1.0]{$U^{s_i}_{k}$}
	\psfrag{p}[][][1.0]{\hspace*{-0.5cm}$New~U^{s_2}_{k}$}
	\psfrag{q}[][][1.0]{\hspace*{-0.5cm}$Old~U^{s_2}_{k}$}
	\psfrag{r}[][][1.0]{$D^{s_2}_k=1$}    
	\subfloat[\label{fig:POT_move_a}]{
		\includegraphics[width=0.8\columnwidth]{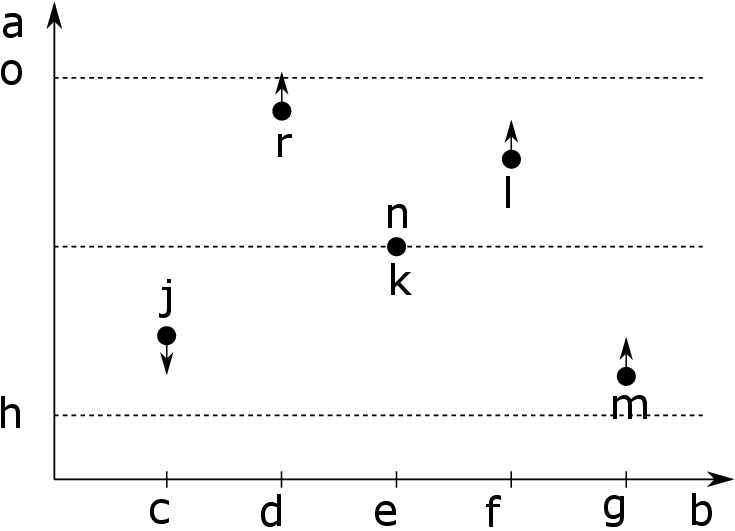}
    }

	\psfrag{a}[][][1.0]{$V^{s_i}_{k+1}$} 
	\psfrag{b}[][][1.0]{$~~s_i$} 
	\psfrag{c}[][][1.0]{$~~s_1$}
	\psfrag{d}[][][1.0]{$~~s_2$}
	\psfrag{e}[][][1.0]{$~~s_3$}
	\psfrag{f}[][][1.0]{$~~s_4$}
	\psfrag{g}[][][1.0]{$~~s_5$}
	\psfrag{h}[][][1.0]{$L^{s_i}_{k+1}$}
	\psfrag{i}[][][1.0]{$D^{s_2}_{k+1}=-1$}
	\psfrag{j}[][][1.0]{$D^{s_1}_{k+1}=-1$}
	\psfrag{k}[][][1.0]{$D^{s_3}_{k+1}=0$}
	\psfrag{l}[][][1.0]{$D^{s_4}_{k+1}=1$}
	\psfrag{m}[][][1.0]{$D^{s_5}_{k+1}=1$}
	\psfrag{n}[][][1.0]{$best~s_i$}
	\psfrag{o}[][][1.0]{$U^{s_i}_{k+1}$}
	\psfrag{p}[][][1.0]{\hspace*{-0.5cm}$~~{\rm New}~U^{s_2}_{k+1}$}
	\psfrag{q}[][][1.0]{\hspace*{-0.5cm}${\rm Old}~U^{s_2}_{k}$}
	\psfrag{r}[][][1.0]{$D^{s_2}_{k+1}=1$}  
    \subfloat[\label{fig:POT_move_b}]{
		\includegraphics[width=0.8\columnwidth]{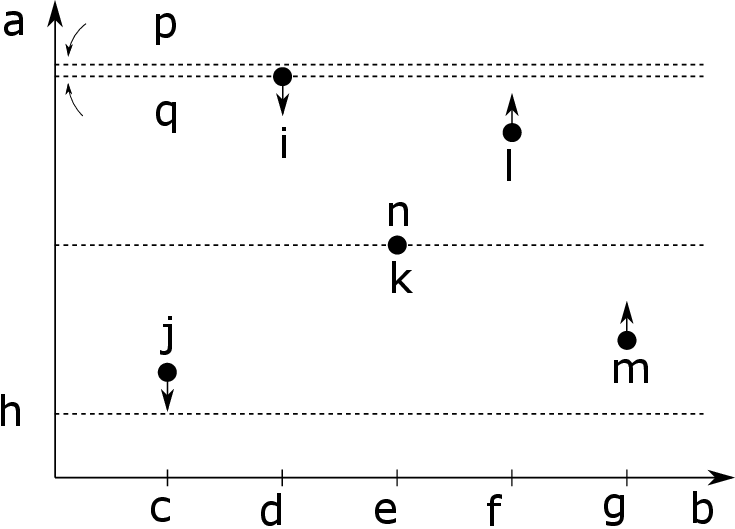}
    }
  	\caption{The Perpetual Orbit technique. (a) The standard behavior; and (b) The increase of $U^{s_2}_{k+1}$ after $V^{s_2}_{k+1}$ reaches the upper bound.}
  \label{fig:POT_move}
\end{figure}

Consider $D^{s_i}_k \in \{-1,1\}$ as the direction of the solution $s_i$ at iteration $k$. If $D^{s_i}_k=1$, $V^{s_i}_k$ has to increase its value by a factor of $(1+\phi)$ unless it reaches its upper bound. In that case, $V^{s_i}_k$ stops increasing its value and $D^{s_i}_k$ becomes $-1$. In a different case, if $D^{s_i}_k=-1$, then  $V^{s_i}_k$ decreases by a factor of $(1-\phi)$ until it reaches its lower bound. At this time,  $V^{s_i}_k$ stops decreasing its value and $D^{s_i}_k$ becomes $1$. In any case, $\phi$ is a small value between $(0.0,0.1]$. The standard behavior of the PO technique  is illustrated in Fig.~\ref{fig:POT_move_a}. Fig.~\ref{fig:POT_move_b} depicts the behavior of the technique when the upper bound is reached.

The pseudo-code of the Perpetual Orbit technique is presented in Algorithm~\ref{alg:POT}.

The PO method can replace the scheduling of the dispersion variable value in ensemble-based algorithms that requires scheduling, which relives the necessity to define the initial value and type of scheduling. In turn, three other parameters are introduced: $\beta$, $\mu$ and $\phi$. However, these parameters do not require scheduling; their values have a well-defined operating range; they do not demand a very tentative tuning process; and they are not dependent on the optimization problem at hand, as shown by the results in this work. Besides, $\beta$ can be eliminated from the algorithm by simple doing it equals 1.

\begin{algorithm}[!ht]
	\caption{Perpetual Orbit technique}
    \label{alg:POT}
    \begin{algorithmic}[1]
    
	\State \textbf{Step 1) Initialization:}\;
	\State \quad - Initialize the chosen algorithm; \;
	\State \quad - Initialize $D^{s_i}_k \in \{-1,1\}$ and $V^{s_i}_k$ randomly,~$\forall s_i \in \Theta$;\;
	\State \quad - Initialize $U^{s_i}_k$ and $L^{s_i}_k$,~$\forall s_i \in \Theta$, according to (\ref{eq:U_L_initial_bounds}); \;

    \item[]

	\State \textbf{Step 2) Generation:}\;
	\State \quad - Generate and assess probing solutions according to the chosen algorithm;\;
	
    \item[]
	
	\State \textbf{Step 3) Evaluation:}\;
	\State \quad - Accept the probing solutions if it follows the chosen algorithm's acceptance criterion;\;
	
	\State \quad - If there is a new best solution, then\;
	\State \quad \quad\quad - Update the Bounds of all solutions according to~(\ref{eq:U_L_initial_bounds});\;
	
    \item[]	
	
	\State \textbf{Step 4) Update:}\;
	\State \quad\textit{4.1) Update the chosen algorithm:}\;
	\State \quad\quad\quad - Update variables of the chosen algorithm;\;

	\State \quad\textit{4.2) Update the Dispersion Variable:}\;

	\State \quad\quad\quad - If $(s_i \neq best( V^{s_i}_k )_{s_i \in \Theta})$ then \;
	\State \quad\quad\quad\quad\quad - If $(D^{s_i}_k > 0)$  then $V^{s_i}_{k+1}$ = $V^{s_i}_k (1+\phi)$;\;
	\State \quad\quad\quad\quad\quad - If $(D^{s_i}_k < 0)$  then $V^{s_i}_{k+1}$ = $V^{s_i}_k (1-\phi)$;\;
	\State \quad\quad\quad ~~Else \;
	\State \quad\quad\quad\quad\quad - $V^{s_i}_{k+1} = V^{s_i}_k$ \;

	\State \quad\textit{4.3) Update the Upper and Lower Bound:}\;
	\State \quad\quad\quad - If $(D^{s_i}_k > 0)$ and $(V^{s_i}_k \geq U^{s_i}_k)$ then \;
		\State \quad\quad \quad\quad\quad - $V^{s_i}_{k+1} = V^{s_i}_{k}$; \quad\quad $D^{s_i}_{k+1} = -1$; \quad\quad $U^{s_i}_{k+1} = U^{s_i}_k(1+\mu);$\;
	
	\State \quad\quad\quad - If $(D^{s_i}_k < 0)$ and $(V^{s_i}_k \leq L^{s_i}_k)$ then \;
		\State \quad\quad \quad\quad\quad - $V^{s_i}_{k+1} = V^{s_i}_{k}$; \quad\quad $D^{s_i}_{k+1} = ~~1$; \quad\quad $L^{s_i}_{k+1} = L^{s_i}_k(1-\mu).$\;

    \item[]
    
	\State \textbf{Step 5) Stopping Criterion:}\;
	\State \quad - Stop if the stopping criterion is met; \;
	\State \quad - Otherwise, go to Step 2. \;
	
    \end{algorithmic}
\end{algorithm}

\section{Perpetual Orbit Coupled Simulated Annealing}
\label{sec:pot_csa}
In CSA, although the control of the variance of the acceptance probabilities makes the optimization more robust with respect to the initial acceptance temperature value, the generation temperature still requires the tuning of a initial value and a monotonic schedule, which is very dependent of the problem at hand. The CSA will stagnate the search for better solutions when its generation temperature is sufficiently diminished, which can lead to convergence to local optima. In order to avoid this, as the CSA is a ensemble-based algorithm and its generation temperature is analogous to a the dispersion variable, the PO technique becomes an interesting strategy. In this case, we use $V^{s_i}_{k+1}$ as $T^{\rm gen_i}_{k+1}$, $D^{s_i}_{k+1}$ as $D^{i}_{k+1}$, $U^{s_i}_{k+1}$ as $U^{i}_{k+1}$ and $L^{s_i}_{k+1}$ as $L^{i}_{k+1}$, i.e., each optimizer has its own generation temperature, direction, upper and lower bounds. When combined to the variance control, it turns the CSA into a parameter-free algorithm, and avoids the stagnation of the optimization in local optima.

The CSA acceptance mechanism was slightly adapted to better work with the PO technique. Although the classical acceptance would work for many cases, it can face problems when the objective function is very smooth or when the initial value of the generation temperature is too small. In both cases, there is a tendency for the generation temperature to be constantly reduced. Low values of the generation temperature generate better solutions whose energy is very close to the current solution, i.e. the technique continues to improve the solution with negligible gains. The less the gain, the less the value of the generation temperature and, consequently, the technique stagnates. This behavior happens because the generation temperature does not have a well-defined scope and any change in its value can have enormous impacts on the optimization process. The solution for these problems is to accept a probing solution only if there is a minimum percentage gain of $\delta$, a small value between $(0.0,0.05]$. This strategy prevents the acceptance of probing solutions with negligible gain, forcing the technique to increase the value of $T^{\rm gen_i}_k$ to an acceptable level where the gains are more expressive. This strategy may not be necessary for other algorithms where the dispersion variable does not need to follow an increasing or decreasing schedule. The pseudo-code for the PO-CSA algorithm is presented in Algorithms~\ref{alg:PO-CSA1} and~\ref{alg:PO-CSA2}. 

\begin{algorithm}[!ht]
	\caption{Part 1 - Perpetual Orbit Coupled Simulated Annealing}
        \label{alg:PO-CSA1}
        
        \begin{algorithmic}[1]
	

	\State \textbf{Step 1) Initialization:}\;
	\State \quad \textit{1.1) CSA}\;
	\State \quad\quad\quad - Assign $m$ random initial solutions to $\Theta$; \;
	\State \quad\quad\quad - Assess the Energy $E(x_i)$, $\forall x_i \in \Theta$; \;
	\State \quad\quad\quad - Set the number of iterations $k=0$;\;
	\State \quad\quad\quad - Set the initial acceptance temperature;\;
	\State \quad\quad\quad - Assign random value to $\alpha$; \;
	\State \quad\quad\quad - Evaluate the coupling term $\gamma$; \;
	\State \quad\quad\quad - Set $\sigma_D^2 = 0.99\frac{(m-1)}{m^2}$.\; 
	\State \quad \textit{1.2) PO}\;
	\State \quad\quad\quad - Assign random values to $\beta$, $\delta$, $\mu$ and $\phi$; \;
	\State \quad\quad\quad - Initialize $D^{i}_k \in \{-1,1\}$ and $T^{\rm gen_i}_k$ randomly,~$\forall i \in \Theta$;\;
	\State \quad\quad\quad - Initialize the $U^{i}_k$ and $L^{i}_k$ according to (\ref{eq:U_L_initial_bounds}); \;

	\item[]
	
	\State \textbf{Step 2) Generation:}\;
	\State \quad\textit{2.1) CSA}\;
	\State \quad\quad\quad - Generate a probing solution $y_i=x_i+\epsilon*T^{\rm tgen_i}_{k}$, $\forall i=1,\cdots,m$,  according to (\ref{eq:cauchy_d}); \;
	\State \quad\quad\quad - Assess the energies $E(y_i), \forall i = 1, ~ \ldots ~ , m$.\;
	
	\item[]
	
	\State \textbf{Step 3) Acceptance:}\;
	\State \quad\textit{3.1) CSA}\;	
	\State \quad\quad\quad - $\forall i=1,\cdots,m$, do $x_i \gets y_i$ if: \;
	\State \quad\quad\quad\quad\quad A) ${y_i}$ is better than $x_i$ in, at the very least, $\delta$; OR \;
	\State \quad\quad\quad\quad\quad B) $A_\Theta > r$, where $r$ is a random variable that was sampled from a uniform distribution $[0,1]$. \;
	\State \quad\textit{3.2) PO}\;	
	\State \quad\quad\quad - If there is a new best solution, then\;
	\State \quad\quad\quad\quad\quad - Update the Bounds of all solutions according to~(\ref{eq:U_L_initial_bounds}).\;
	
    \end{algorithmic}
\end{algorithm}

\begin{algorithm}[!ht]
	\caption{Part 2 - Perpetual Orbit Coupled Simulated Annealing}
        \label{alg:PO-CSA2}
        \begin{algorithmic}[1]
	
    
	\State \textbf{Step 4) Update:}\;
	\State \quad\textit{4.1) CSA}\;
	\State \quad\quad\quad - Calculate $\sigma^2$; \;
	\State \quad\quad\quad - Adjust $T^{\rm ac}_k$ according to the following rules: \;
	\State \quad\quad\quad\quad \quad- If $\sigma^2 < \sigma_D^2$ then $T^{\rm ac}_{k+1} = T^{\rm ac}_k(1-\alpha)$; \;
	\State \quad\quad\quad\quad\quad - If $\sigma^2 > \sigma_D^2$ then $T^{\rm ac}_{k+1} = T^{\rm ac}_k(1+\alpha)$. \;
        \State \quad\quad\quad - Evaluate the coupling term $\gamma$; \;
   	\State \quad\textit{4.2) PO}\;
   	\State \quad\quad\quad\textit{4.2.1) Update the Generation Temperature:}\;
	\State \quad\quad\quad\quad\quad - If $(i \neq best( T^{\rm gen_i}_k )_{i \in \Theta})$ then \;
	\State \quad\quad\quad\quad\quad\quad\quad - If $(D^{i}_k > 0)$  then $T^{\rm gen_i}_{k+1}$ = $T^{\rm gen_i}_k (1+\phi)$;\;
	\State \quad\quad\quad\quad\quad\quad\quad - If $(D^{i}_k < 0)$  then $T^{\rm gen_i}_{k+1}$ = $T^{\rm gen_i}_k (1-\phi)$;\;
	\State \quad\quad\quad\quad\quad ~~Else \;
	\State \quad\quad\quad\quad\quad\quad\quad - $T^{\rm gen_i}_{k+1} = T^{\rm gen_i}_k$ \;

	\State \quad\quad\quad\textit{4.2.2) Update the Upper and Lower Bound:}\;
	\State \quad\quad\quad\quad\quad - If $(D^{i}_k > 0)$ and $(T^{\rm gen_i}_k \geq U^{i}_k)$ then \;
		\State \quad\quad\quad\quad\quad\quad\quad - $T^{\rm gen_i}_{k+1} = T^{\rm gen_i}_{k}$; \quad $D^{i}_{k+1} = -1$; \;
		\State \quad\quad\quad\quad\quad\quad\quad - $U^{i}_{k+1} = U^{i}_k(1+\mu);$\;
	
	\State \quad\quad\quad\quad\quad - If $(D^{i}_k < 0)$ and $(T^{\rm gen_i}_k \leq L^{i}_k)$ then \;
		\State \quad\quad\quad\quad\quad\quad\quad - $T^{\rm gen_i}_{k+1} = T^{\rm gen_i}_{k}$; \quad $D^{i}_{k+1} = +1$; \;
		\State \quad\quad\quad\quad\quad\quad\quad - $L^{i}_{k+1} = L^{i}_k(1-\mu);$\;
	\State \quad\quad\quad\quad\quad - Increment $k$ by 1.\;

        \item[]
    
	\State \textbf{Step 5) Stopping Criterion:}\;
	\State \quad - Stop if the stopping criterion is met; \;
	\State \quad - Otherwise, go to Step 2. \;
	
    \end{algorithmic}
\end{algorithm}

The PO-CSA algorithm starts with \textit{Initialization} step by setting initial variables, such as the initial solutions,  the generation and acceptance temperatures, the upper and lower bounds of $T^{\rm gen_i}_k$, and the direction of each solution. The energy of each initial solution is assessed and the coupling term is calculated.

In \textit{Generation}, the PO-CSA behaves like the CSA, by generating probing solutions and assessing their energies. This generation needs a statistical distribution, which in turn needs a dispersion variable. In the CSA and PO-CSA, the generation temperature is used as the dispersion variable. 

In \textit{Acceptance}, a better solution is always accepted and a worse probing solution is accepted if it is better than the current solution by a percentage gain of $\delta$, at least, or, with probability $A_\Theta(\gamma, x_i \to y_i)$. If any accepted solution is the new best solution, each optimizer has to update its upper and lower bounds $U_k^{i}$ and $L_k^{i}$, respectively, according to (\ref{eq:U_L_initial_bounds}).

Then, the algorithm need to \textit{Update} the acceptance temperature and the generation temperatures. The first one is updated by following the rule described in~(\ref{eq:cond_var}). The coupling term is reevaluated. The generation temperatures of each solution, $T^{\rm gen_i}_k$, follow the perpetual orbit, where they move around the best solution's generation temperature. The upper or lower bound is updated if a solution's generation temperature reaches it. Any optimizer that finds the best solutions keeps its generation temperature unchanged and becomes the new reference for all other generation temperatures. Each optimizer, then, has to update its upper and lower bounds $U_k^{i}$ and $L_k^{i}$, respectively.

The algorithm finishes its execution when the \textit{Stopping Criterion} is met. Otherwise, the loop restarts through a new execution of the \textit{Generation} step.

\section{Experiments \& Results}
\label{sec:exp_results}

We performed several experiments to assess the optimization potential of the proposed PO-CSA algorithm. We performed the experiments on a cluster with 64 nodes, each one with 2 Intel Xeon Sixteen-Core E5-2698v3 CPUs at 2.3 GHz, and 128 GB of RAM DDR4, from the High Performance Computing Center at UFRN (NPAD/UFRN). We used the same 14 objective functions used in the experiments with the original CSA, better described in Appendix~\ref{appendix}. The initialization of the generation temperatures and its schedule as well as the method for generating probing solutions are explained in the following subsection. The results of the experiments are presented in section~\ref{subsec:results}.

\subsection{General Settings}
\label{subsec:general_settings}
In the steps 2 and 4 of the CSA, it is necessary to generate and assess a new probing solutions, and decrease a generation temperature $T^{\rm gen}_k$, previously initialized with a value $T^{\rm gen}_0$, according to a chosen monotonic schedule, respectively. SA-based algorithms tend to be very sensitive to different schedules and initial parameters. The impact of these choices goes beyond the convergence of the algorithm, directly impacting the quality of the final solution. In the last two of the four sets of experiments performed in this work, two initialization procedures were used for the CSA to compare the resulting behavior with the PO-CSA: R-CSA and B-CSA. They work exactly as the original CSA, however the first one had its initial generation temperature initialized as a random number between $\left[0,100\right]$, just like it was done for the PO-CSA, while the B-CSA shows the best results from the best initial generation temperature chosen from a predefined set of initial generation temperatures
\begin{equation}
\label{eq::init_tgen_bcsa}
		T^{\rm gen}_{0} \in \{0.001,0.01,0.1,1,10,100,1000\},
\end{equation}
i.e. for each of the experiments that is presented here, the B-CSA shows the best solution among the solutions found by the CSA with these seven different initial generation temperatures.

The R-CSA and the B-CSA were subjected to the same generation temperature schedule~\cite{szu1987fast}:
\begin{equation}
		T^{\rm gen}_{k+1} = \frac{T^{\rm gen}_0}{k+1}.
\end{equation}
By choosing the same generation temperature schedule to both algorithms, we hope to better demonstrate the behavior of the PO-CSA algorithm.  There was no study to conclude which schedule would be the best suit to run each function on the B-CSA and the R-CSA.

In the R-CSA, B-CSA and PO-CSA, we set the increase/decrease acceptance temperature factor $\alpha=0.05$ and the desired variance of the acceptance probabilities to $99\%$ of its maximum value. Exclusively for the PO-CSA, we used the boundary multiplier $\beta=10$, the increase/decrease boundary factor $\mu=\alpha$ and percentage gain $\delta=0.001$ ($0.1\%$) for all functions and experiments.

\subsection{Results}
\label{subsec:results}
The experiments were divided in four sets.  The first set compares the quality of solution of the PO-CSA and four reference algorithms. The second set shows the behavior of the generation temperature when the PO works combined with the CSA. The third set demonstrates the robustness of the quality of the final solutions of the PO-CSA comparing it to the CSA with two different initialization procedures: R-CSA and B-CSA. The fourth set shows the resilience of the proposed method for a larger number of function evaluations by comparing it to the R-CSA and the B-CSA cases using the function where PO-CSA performed worst in the third set of experiments. These experiments and their results are detailed as follows.

\subsubsection{Comparing PO-CSA to Classical Algorithms}
\label{subsubsec:comparation}
In this set of experiment, we compare theperformance of the PO-CSA to the average solutions found by four of the most used heuristics in scientific applications~\cite{civicioglu2013conceptual}: the Cuckoo Search via L\'{e}vy Flight (CS)~\cite{yang2009cuckoo,galvez2014cuckoo}, Differential Evolution (DE)~\cite{storn1997differential,song2023dynamic,piotrowski2023particle}, Particle Swarm Optimization (PSO)~\cite{kenny1995particle,piotrowski2023particle}, and Genetic Algorithm (GA)~\cite{kuhlemann2020genetic,wright1991genetic,holland1992adaptation,eshelman1993real}. The stopping criterion was $D\times10^6$ function evaluations, where $D$ is the dimension of the problem. The population size for CS, DE and PSO was 50~\cite{civicioglu2013conceptual}. As in~\cite{sarmady2007investigation} the recommended population size of the GA is no greater than 100, and in~\cite{civicioglu2013conceptual} the recommended the population size is 50 for CS, DE, and PSO, we decided to use the same population size for GA.  The number of optimizers used in the the PO-CSA was $D$. The results are averaged across 25 runs. 

In order to avoid an incessant search for better parameters for the classic algorithms, we have used the following values of algorithmic control parameters according to the references:
\begin{itemize}
\item CS: $\beta=1.5$ and $p_0=0.25$~\cite{yang2010engineering,civicioglu2013conceptual};
\item DE: $DE/rand/1$ mutation strategy with $F=0.50$ and $Cr=0.90$~\cite{neri2010recent,das2011differential,karaboga2009comparative,civicioglu2013conceptual};
\item PSO: $C_1=C_2=1.80$ and $\omega=0.60$~\cite{karaboga2009comparative,civicioglu2013conceptual};
\item GA: Blend Crossover (BLX-$\alpha$) with $\alpha=0.5$~\cite{herrera2003taxonomy,takahashi2001crossover};
\item PO-CSA: $\beta=10$ , $\delta=0.001$ ($0.1\%$) and $\mu=\alpha=0.05$.
\end{itemize}

The average results are shown in Table~\ref{tab:quality_solution_comparison}. The average solutions found by the PO-CSA are equal or better than the other algorithms in 85.71\% of the cases, against 2.38\% for the CS and GA, 9.52\% for the DE and 17.85\% for the PSO. There is a tie of, at least, 2 metaheuristics in 15.47\% of the cases.

\begin{sidewaystable}
\caption{Statistical average solution of the CS, DE, GA, PSO and PO-CSA for functions $f_1-f_{14}$ with $D \in \{5, 10, 15, 20, 25, 30\}$ and the number of optimizer equals $D$. The stopping criterion is  $1D\times10^6$ function evaluations. The population size for Cuckoo, DE, GA, PSO is 50. The number of optimizers of the PO-CSA is $D$. The results are averaged across 25 runs. The best results for each dimension in a given number of function evaluations are highlighted in bold.}\label{tab:quality_solution_comparison}
\tiny
\begin{tabular}{c|c|c|c|c|c|c||c|c|c|c|c|c|c}
	\hline
		\multirow{2}{*}{Function} & \multirow{2}{*}{Dimension} &  \multicolumn{5}{c||}{Algorithm} & \multirow{2}{*}{Function} & \multirow{2}{*}{Dimension} &  \multicolumn{5}{c}{Algorithm}	\\ 
        \cmidrule{3-7}
        \cmidrule{10-14}
            & & CS & DE & GA & PSO & PO-CSA & & &  CS & DE & GA & PSO & PO-CSA \\
  	    \hline
    	\multirow{6}{*}{$f_1$} & $D=05$ & 1.21E-06 & \textbf{0.00E+00} & 2.91E-05 & \textbf{0.00E+00} & \textbf{0.00E+00} &	\multirow{6}{*}{$f_8$} & $D=05$ & \textbf{8.56E-02} & 7.14E+01 & 8.60E-02 & 9.96E+01 & \textbf{8.56E-02}\\
    	& $D=10$ & 9.47E-05 & 2.20E-04 & 1.43E-04 & \textbf{0.00E+00} & \textbf{0.00E+00} & & $D=10$& \textbf{1.71E-01} & 3.78E+02 & 1.72E-01 & 4.85E+02 & \textbf{1.71E-01}\\
    	& $D=15$ & 3.64E-03 & 3.60E-01 & 5.05E-04 & \textbf{0.00E+00} & \textbf{0.00E+00} & & $D=15$& 2.67E-01 & 1.05E+03 & 2.58E-01 & 1.13E+03 & \textbf{2.57E-01}\\
    	& $D=20$ & 1.12E-01 & 1.03E+01 & 5.17E-04 & \textbf{0.00E+00} & \textbf{0.00E+00} & & $D=20$& 6.22E-01 & 1.74E+03 & \textbf{3.43E-01} & 2.18E+03 & \textbf{3.43E-01}\\
    	& $D=25$ & 2.50E+00 & 1.74E+00 & 2.71E-04 & \textbf{0.00E+00} & \textbf{0.00E+00} & & $D=25$& 7.93E+00 & 2.48E+03 & \textbf{4.28E-01} & 3.07E+03 & \textbf{4.28E-01}\\
    	& $D=30$ & 4.14E+01 & 2.70E+00 & 4.93E-04 & \textbf{0.00E+00} & \textbf{0.00E+00} & & $D=30$& 1.37E+02 & 3.47E+03 & 5.15E-01 & 4.46E+03 & \textbf{5.14E-01}\\
    	\hline
    	\multirow{6}{*}{$f_2$} & $D=05$ 	& 4.07E-02 & 8.33E-01 & 1.48E+00 & \textbf{1.70E-28} & 3.88E-19 & \multirow{6}{*}{$f_9$} & $D=05$ & 1.61E+00 & 2.90E-07 & 2.13E+00 & 3.56E-01 & \textbf{4.44E-16}\\
    	& $D=10$ & 1.67E+00 & 7.13E+00 & 5.80E+00 & 3.19E-01 & \textbf{8.86E-10} & & $D=10$& 2.79E+00 & 3.96E-01 & 3.21E+00 & 7.25E-01 & \textbf{4.00E-15}\\
    	& $D=15$ & 8.67E+00 & 1.26E+01 & 9.27E+00 & 6.39E-01 & \textbf{2.25E-05} & & $D=15$& 3.26E+00 & 1.05E+00 & 3.27E+00 & 1.67E+00 & \textbf{1.64E-13}\\ 
    	& $D=20$ & 2.59E+01 & 1.78E+01 & 1.44E+01 & 3.19E-01 & \textbf{1.47E-03} & & $D=20$& 3.78E+00 & 2.21E+00 & 3.67E+00 & 1.23E+00 & \textbf{7.23E-02}\\
    	& $D=25$ & 3.79E+01 & 2.26E+01 & 2.42E+01 & 6.38E-01 & \textbf{5.89E-02} & & $D=25$& 3.86E+00 & 2.92E+00 & 3.99E+00 & \textbf{1.03E+00} & 1.16E+00\\
    	& $D=30$ & 7.67E+01 & 2.71E+01 & 2.69E+01 & 4.78E-01 & \textbf{2.16E-01} & & $D=30$& 4.39E+00 & 3.20E+00 & 4.35E+00 & \textbf{1.48E+00} & 1.78E+00\\ 
    	\hline
    	\multirow{6}{*}{$f_3$} & $D=05$ 	& 6.80E-04 & \textbf{4.44E-16} & 2.60E-03 & \textbf{4.44E-16} & \textbf{4.44E-16} &	\multirow{6}{*}{$f_{10}$} & $D=05$ & 3.44E-01 & 1.79E-01 & 5.92E-01 & 1.14E-01 & \textbf{1.73E-02}\\
    	& $D=10$ & 4.22E-03 & 2.51E-01 & 4.75E-03 & \textbf{3.43E-15} & 4.00E-15 & & $D=10$& 5.63E-01 & 1.37E-01 & 6.46E-01 & 1.75E-01 & \textbf{3.40E-02}\\
    	& $D=15$ & 4.31E-02 & 1.42E+00 & 4.32E-03 & 4.14E-15 & \textbf{4.00E-15} & & $D=15$& 3.93E-01 & 1.27E-01 & 3.47E-01 & 1.24E-01 & \textbf{1.40E-02}\\
    	& $D=20$ & 1.55E+00 & 2.68E+00 & 4.18E-03 & \textbf{5.28E-15} & 7.55E-15 & & $D=20$& 6.98E-01 & 2.36E-01 & 4.26E-01 & 5.12E-02 & \textbf{1.17E-02}\\
    	& $D=25$ & 2.70E+00 & 3.35E+00 & 6.35E-03 & 1.85E-01 & \textbf{7.55E-15} & & $D=25$& 1.08E+00 & 3.23E-01 & 4.06E-01 & 3.43E-02 & \textbf{9.51E-03}\\
    	& $D=30$ & 4.31E+00 & 4.09E+00 & 5.48E-03 & \textbf{7.41E-15} & 7.55E-15 & & $D=30$& 1.42E+00 & 1.98E-01 & 3.46E-01 & 1.40E-02 & \textbf{3.41E-03}\\
    	\hline
    	\multirow{6}{*}{$f_4$} & $D=05$ 	& 8.98E-02 & 6.41E-02 & 1.47E-02 & 1.91E-02 & \textbf{2.27E-04}  &	\multirow{6}{*}{$f_{11}$} & $D=05$ & 5.13E-01 & \textbf{3.61E-05} & 6.23E-01 & 2.46E-01 & 4.77E-05\\
    	& $D=10$ & 2.10E-01 & 1.21E-01 & 5.00E-02 & 6.38E-02 & \textbf{7.40E-03} & & $D=10$& 2.90E+00 & \textbf{1.61E-01} & 1.98E+00 & 7.37E-01 & 1.86E-01\\
    	& $D=15$ & 9.12E-02 & 6.79E-02 & 6.03E-02 & 4.05E-02 & \textbf{0.00E+00} & & $D=15$& 4.76E+00 & \textbf{7.69E-01} & 4.02E+00 & 1.30E+00 & 1.15E+00\\ 
    	& $D=20$ & 3.60E-01 & 1.10E-01 & 5.09E-02 & 2.32E-02 & \textbf{0.00E+00} & & $D=20$& 7.08E+00 & \textbf{2.21E+00} & 5.25E+00 & 2.71E+00 & 2.50E+00\\
    	& $D=25$ & 1.00E+00 & 1.48E-01 & 6.36E-02 & 9.93E-03 & \textbf{0.00E+00} & & $D=25$& 8.48E+00 & \textbf{2.86E+00} & 8.03E+00 & 3.27E+00 & 4.41E+00\\
    	& $D=30$ & 1.43E+00 & 2.40E-01 & 6.50E-02 & 1.12E-02 & \textbf{0.00E+00} & & $D=30$& 1.31E+01 & \textbf{5.16E+00} & 9.50E+00 & 6.11E+00 & 7.27E+00\\
    	\hline
    	\multirow{6}{*}{$f_5$} & $D=05$ 	& 3.11E-02 & \textbf{0.00E+00} & 3.19E-03 & \textbf{0.00E+00} & \textbf{0.00E+00} &	\multirow{6}{*}{$f_{12}$} & $D=05$ & 2.21E+00 & 3.79E+00 & 4.53E+00 & 1.79E+00 & \textbf{6.21E-08}\\
    	& $D=10$ & 2.66E-01 & 4.79E-03 & 7.04E-03 & \textbf{0.00E+00} & \textbf{0.00E+00} & & $D=10$& 1.48E+01 & 1.27E+01 & 1.46E+01 & 8.52E+00 & \textbf{2.85E+00}\\
    	& $D=15$ & 7.65E-01 & 2.66E-01 & 1.79E-02 & 5.18E-04 & \textbf{0.00E+00} & & $D=15$& 3.08E+01 & 2.25E+01 & 2.83E+01 & 1.69E+01 & \textbf{8.05E+00}\\
    	& $D=20$ & 1.81E+00 & 1.01E+00 & 2.88E-02 & 1.26E-01 & \textbf{0.00E+00} & & $D=20$& 4.49E+01 & 2.93E+01 & 3.80E+01 & 2.09E+01 & \textbf{1.50E+01}\\
    	& $D=25$ & 3.84E+00 & 7.89E+00 & 3.33E-02 & 2.10E-01 & \textbf{0.00E+00} & & $D=25$& 5.30E+01 & 3.99E+01 & 4.26E+01 & 3.46E+01 & \textbf{2.02E+01}\\
    	& $D=30$ & 6.89E+00 & 2.74E+00 & 4.64E-02 & 6.53E-01 & \textbf{0.00E+00} & & $D=30$& 7.38E+01 & 4.52E+01 & 4.93E+01 & 3.98E+01 & \textbf{3.19E+01}\\
    	\hline
    	\multirow{6}{*}{$f_6$} & $D=05$ 	& 8.14E-07 & 7.48E-01 & 1.55E-05 & 4.78E-01 & \textbf{0.00E+00 }&	\multirow{6}{*}{$f_{13}$} & $D=05$ & 1.83E+00 & 2.75E+00 & 4.16E+00 & 1.88E+00 & \textbf{1.32E-06}\\
    	& $D=10$ & 4.90E-05 & 5.59E+00 & 6.01E-05 & 4.06E+00 & \textbf{0.00E+00} & & $D=10$& 7.00E+00 & 1.11E+01 & 1.87E+01 & 7.24E+00 & \textbf{2.47E+00}\\
    	& $D=15$ & 2.83E-03 & 1.23E+01 & 1.76E-04 & 9.91E+00 & \textbf{0.00E+00} & & $D=15$& 1.80E+01 & 1.52E+01 & 3.30E+01 & 1.45E+01 & \textbf{5.00E+00}\\
    	& $D=20$ & 8.15E-01 & 2.19E+01 & 2.37E-04 & 1.72E+01 & \textbf{0.00E+00} & & $D=20$& 2.94E+01 & 2.45E+01 & 4.39E+01 & 2.64E+01 & \textbf{1.20E+01}\\ 
    	& $D=25$ & 7.12E+00 & 3.17E+01 & 1.66E-04 & 2.69E+01 & \textbf{0.00E+00} & & $D=25$& 4.49E+01 & 3.02E+01 & 5.06E+01 & 3.05E+01 & \textbf{1.50E+01}\\
    	& $D=30$ & 1.88E+01 & 4.31E+01 & 3.74E-04 & 3.62E+01 & \textbf{0.00E+00} & & $D=30$& 6.63E+01 & 3.86E+01 & 5.88E+01 & 4.20E+01 & \textbf{2.10E+01}\\
    	\hline
    	\multirow{6}{*}{$f_7$} & $D=05$ 	& 7.25E-07 & 8.70E-01 & 3.18E-05 & 5.20E-01 & \textbf{0.00E+00} &	\multirow{6}{*}{$f_{14}$} & $D=05$ & 2.17E-01 & 1.20E+02 & 2.94E+02 & 1.36E+02 & \textbf{8.57E-02}\\
    	& $D=10$ & 5.06E-05 & 6.53E+00 & 7.14E-05 & 1.88E+00 & \textbf{0.00E+00} & & $D=10$& 4.77E+02 & 8.89E+02 & 1.18E+03 & 8.96E+02 & \textbf{3.27E+00}\\
    	& $D=15$ & 3.20E-03 & 1.21E+01 & 1.36E-04 & 3.00E+00 & \textbf{0.00E+00} & & $D=15$& 1.37E+03 & 1.89E+03 & 1.92E+03 & 1.74E+03 & \textbf{6.48E+02}\\
    	& $D=20$ & 7.31E-01 & 1.88E+01 & 8.98E-05 & 4.52E+00 & \textbf{0.00E+00} & & $D=20$& 2.10E+03 & 3.04E+03 & 2.84E+03 & 2.71E+03 & \textbf{1.61E+03}\\
    	& $D=25$ & 5.80E+00 & 2.90E+01 & 2.23E-04 & 7.64E+00 & \textbf{0.00E+00} & & $D=25$& 2.44E+03 & 3.47E+03 & 3.26E+03 & 3.52E+03 & \textbf{2.08E+03}\\
    	& $D=30$ & 1.51E+01 & 3.92E+01 & 1.92E-04 & 7.36E+00 & \textbf{0.00E+00} & & $D=30$& 3.60E+03 & 5.03E+03 & 4.17E+03 & 4.57E+03 & \textbf{3.14E+03}\\
    	\hline
    \hline
    \end{tabular}
\end{sidewaystable}

\subsubsection{Analysis of the Generation Temperature Control}
\label{subsubsec:quality_solution_PO}
\begin{figure}[!ht]
	\centering
	\psfragscanon
	\psfrag{x}[][][0.7]{$Iterations$}
	\psfrag{Tgen=0_001}[][][0.5]{$T^{\rm gen}_0=0.001$}
	\psfrag{Tgen=1}[][][0.5]{$T^{\rm gen}_0=1$}
	\psfrag{Tgen=1000}[][][0.5]{$T^{\rm gen}_0=1000$}
	\psfrag{y}[][][0.7]{$Reference~of~the~Generation~temperatures$}
    \subfloat[\label{fig:f3_data_bestTgens}]{
		\includegraphics[width=0.7\columnwidth]{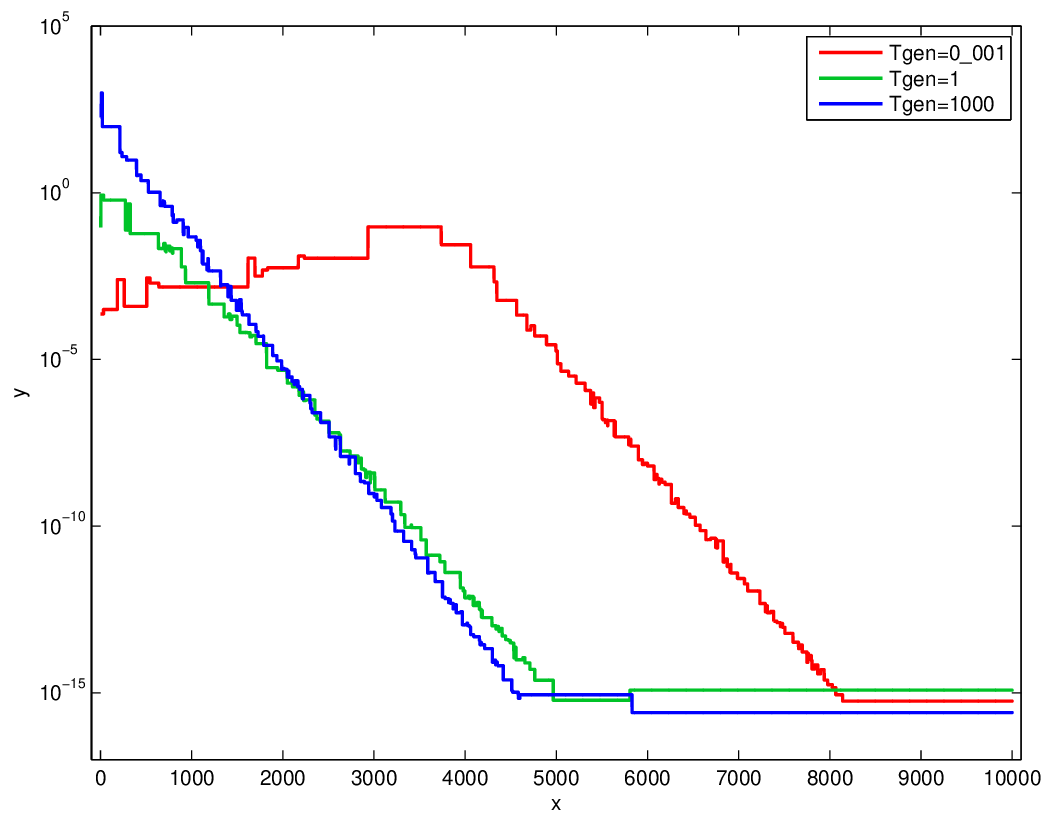}
    }

	\psfrag{y}[][][0.7]{$Generation~temperature$}
    \subfloat[\label{fig:f3_data_allTgens}]{
		\includegraphics[width=0.7\columnwidth]{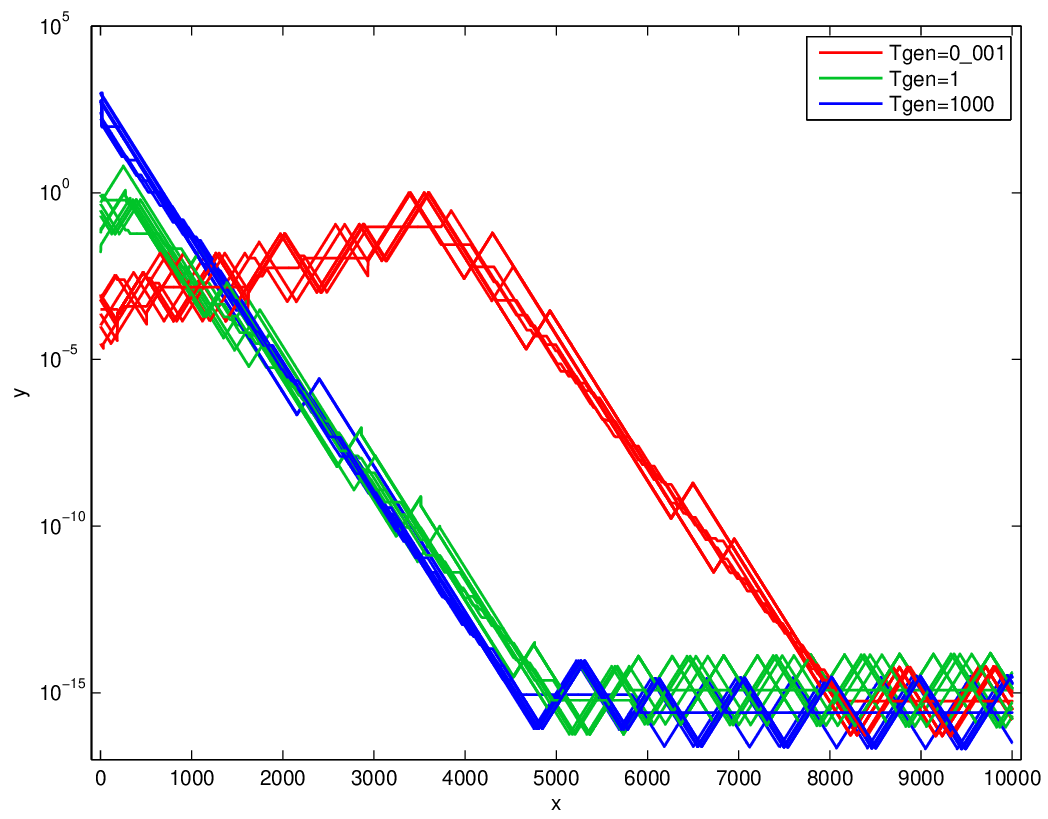}
    }    
  	\caption{The PO-CSA behavior on $f_3$ along 10 thousand iterations with $D=10$ for three different initial generation temperatures. The number of optimizers equals $D$. After a transient, all three runs converge to the same level of energy. (a) The reference for the generation temperature, (b) The generation temperature for all optimizers.}
  \label{fig:f3_data2}
\end{figure}

In this section, we analyze the effects of the generation temperature control through PO on the CSA, as detailed in Section~\ref{sec:pot_csa}. Fig.~\ref{fig:f3_data2} presents the PO-CSA generation temperature's behavior on $f_3$ for three different initial generation temperature: 0.001, 1 and 1000. Fig.~\ref{fig:f3_data_bestTgens} depicts their reference generation temperature, i.e. the value where the generation temperature of the other solutions has to move around it. This reference is the generation temperature of the best solution, as described in section~\ref{sec:pot}. Fig.~\ref{fig:f3_data_allTgens} shows all the generations temperatures around the reference generation temperatures. Note that the convergence of PO-CSA for $T^{\rm gen}_0=0.001$ is slower than the others because its initial value is too small. The control needs a few thousands of iterations before the energy starts descending in a rate similar to the other initial temperatures. The same happens to $T^{\rm gen}_0=1000$, but, in this case, it is because its initial value is too large. Although this behavior reflects in the quality of solution, all three runs converge to the same level of energy and generation temperature after a transient.

\begin{figure}[!ht]
	\centering
	\psfragscanon
	\psfrag{x}[][][0.7]{$Iterations$} 
	\psfrag{y}[][][0.7]{$Best~Energy$} 
	\psfrag{Tgen=0_001}[][][0.5]{$T^{\rm gen}_0=0.001$}
	\psfrag{Tgen=1}[][][0.5]{$T^{\rm gen}_0=1$}
	\psfrag{Tgen=1000}[][][0.5]{$T^{\rm gen}_0=1000$}
	\label{fig:f3_data_AllEnergies}
	\includegraphics[width=0.8\columnwidth]{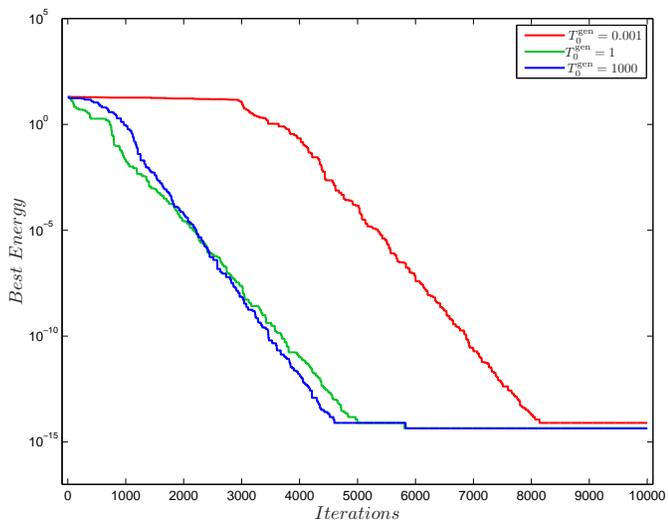}
  	\caption{The best PO-CSA energies on $f_3$ along 10 thousand iterations with $D=10$ for three different initial generation temperatures. The number of optimizers equals $D$. After a short transient, all three runs converge to the same level of energy.}
  \label{fig:f3_data1}
\end{figure}

Fig.~\ref{fig:f3_data1} presents the best energy generated by results of the generation temperature control in Fig.~\ref{fig:f3_data2}. After a transient, all three runs converge to the same level of energy. Observe that the generation temperature convergence impacts on the quality of solutions. As long as the algorithms does not find a good reference generation temperature, the quality of the solution does not improve significantly. However, after reaching this reference, there is a sufficiently great improvement in the quality of the solution.

\subsubsection{Robustness of Initialization}
\label{subsubsec:quality_solution}
In this set of experiments, we compared the quality of the solutions found by the PO-CSA to those found by the R-CSA and B-CSA. The R-CSA has its initial generation temperature as a random number between [0,100], while the B-CSA shows the best results from the best initial generation temperature chosen from a predefined set of initial generation temepratures, as defined in~(\ref{eq::init_tgen_bcsa}). 

n order to perform these experiments, we used six different values for the dimension of the problems $D \in \{5, 10, 15, 20, 25, 30\}$ as the same number of SA optimizers. The stopping criterion for the three algorithms was 1, 2, 4 and 8 million of function evaluations per each optimizer and results were averaged across 25 runs. The mean results are shown in Table~\ref{tab:quality_solution1} and Table~\ref{tab:quality_solution2}.

\begin{sidewaystable}
\caption{Statistical mean solution of the R-CSA, the B-CSA and the PO-CSA for functions $f_1-f_7$ with $D \in \{5, 10, 15, 20, 25, 30\}$ and the number of optimizer equals $D$. The stopping criterion goes from  $1\times10^6$ to $8\times10^6$ function evaluations per optimizer. The results are averaged across 25 runs. The best results for each dimension in a given number of function evaluations are highlighted in bold. }\label{tab:quality_solution1}
\tiny

\begin{tabular}{c|c|cccc|cccc|cccc}
	\hline
		\multicolumn{2}{c|}{Algorithm}	& \multicolumn{4}{c|}{R-CSA} & \multicolumn{4}{c|}{B-CSA} & \multicolumn{4}{c}{PO-CSA} \\
		Function & Dimension & $1\times10^6$ & $2\times10^6$ & $4\times10^6$ & $8\times10^6$	& $1\times10^6$ & $2\times10^6$ & $4\times10^6$ & $8\times10^6$	& $1\times10^6$ & $2\times10^6$ & $4\times10^6$ & $8\times10^6$	\\
  	    \hline
    	\multirow{6}{*}{$f_1$} & $D=05$ & \textbf{0.00E+00} &\textbf{0.00E+00} &\textbf{0.00E+00} &\textbf{0.00E+00} &\textbf{0.00E+00} &\textbf{0.00E+00} &\textbf{0.00E+00} &\textbf{0.00E+00} &\textbf{0.00E+00} &\textbf{0.00E+00} &\textbf{0.00E+00} &\textbf{0.00E+00}  \\
    	& $D=10$ & \textbf{0.00E+00} &\textbf{0.00E+00} &\textbf{0.00E+00} &\textbf{0.00E+00} &\textbf{0.00E+00} &\textbf{0.00E+00} &\textbf{0.00E+00} &\textbf{0.00E+00} &\textbf{0.00E+00} &\textbf{0.00E+00} &\textbf{0.00E+00} &\textbf{0.00E+00}  \\
    	& $D=15$ & \textbf{0.00E+00} &\textbf{0.00E+00} &\textbf{0.00E+00} &\textbf{0.00E+00} &\textbf{0.00E+00} &\textbf{0.00E+00} &\textbf{0.00E+00} &\textbf{0.00E+00} &\textbf{0.00E+00} &\textbf{0.00E+00} &\textbf{0.00E+00} &\textbf{0.00E+00}  \\
    	& $D=20$ & \textbf{0.00E+00} &\textbf{0.00E+00} &\textbf{0.00E+00} &\textbf{0.00E+00} &\textbf{0.00E+00} &\textbf{0.00E+00} &\textbf{0.00E+00} &\textbf{0.00E+00} &\textbf{0.00E+00} &\textbf{0.00E+00} &\textbf{0.00E+00} &\textbf{0.00E+00}  \\
    	& $D=25$ & \textbf{0.00E+00} &\textbf{0.00E+00} &\textbf{0.00E+00} &\textbf{0.00E+00} &\textbf{0.00E+00} &\textbf{0.00E+00} &\textbf{0.00E+00} &\textbf{0.00E+00} &\textbf{0.00E+00} &\textbf{0.00E+00} &\textbf{0.00E+00} &\textbf{0.00E+00}  \\
    	& $D=30$ & \textbf{0.00E+00} &\textbf{0.00E+00} &\textbf{0.00E+00} &\textbf{0.00E+00} &\textbf{0.00E+00} &\textbf{0.00E+00} &\textbf{0.00E+00} &\textbf{0.00E+00} &\textbf{0.00E+00} &\textbf{0.00E+00} &\textbf{0.00E+00} &\textbf{0.00E+00}  \\
    	\hline
\multirow{6}{*}{$f_2$} & $D=05$ 	& 7.57E-06 &1.89E-06 &3.00E-07 &6.65E-08 &3.08E-12 &5.83E-13 &7.99E-14 &1.91E-14 &\textbf{3.88E-19} &\textbf{3.62E-22} &\textbf{1.80E-28} &\textbf{0.00E+00} \\
    	& $D=10$ & 6.13E-04 &8.07E-05 &2.13E-05 &5.58E-06 &3.41E-09 &3.03E-10 &2.97E-11 &8.45E-13 &\textbf{8.86E-10} &\textbf{6.82E-15} &\textbf{6.51E-24} &\textbf{2.49E-29} \\
    	& $D=15$ & 5.11E-03 &7.20E-04 &1.07E-04 &1.81E-05 &\textbf{3.51E-07} &2.01E-08 &2.96E-10 &6.24E-11 &2.25E-05 &\textbf{8.03E-09} &\textbf{2.66E-14} &\textbf{3.69E-22} \\
    	& $D=20$ & 8.07E-03 &2.48E-03 &4.76E-04 &1.38E-04 &\textbf{1.71E-06} &\textbf{1.02E-07} &\textbf{7.72E-10} &1.05E-10 &1.47E-03 &2.13E-05 &3.12E-09 &\textbf{5.53E-15} \\
    	& $D=25$ & 3.64E-02 &9.47E-03 &1.24E-03 &1.38E-04 &\textbf{5.01E-06} &\textbf{4.22E-07} &\textbf{1.82E-09} &\textbf{4.00E-10} &5.89E-02 &5.47E-02 &4.82E-04 &5.54E-09 \\
    	& $D=30$ & 8.07E-02 &1.29E-02 &2.22E-03 &6.23E-04 &\textbf{1.11E-05} &\textbf{9.28E-07} &\textbf{1.03E-07} &\textbf{8.92E-09} &2.16E-01 &1.19E-01 &3.66E-03 &1.09E-06 \\
    	\hline
    	\multirow{6}{*}{$f_3$} & $D=05$ & \textbf{4.44E-16} & \textbf{4.44E-16} & \textbf{4.44E-16} & \textbf{4.44E-16} & \textbf{4.44E-16} & \textbf{4.44E-16} & \textbf{4.44E-16} & \textbf{4.44E-16} & \textbf{4.44E-16} & \textbf{4.44E-16} & \textbf{4.44E-16} & \textbf{4.44E-16} \\
    	& $D=10$ & 4.00E-15 & 4.00E-15 & 3.85E-15 & 3.71E-15 & \textbf{3.85E-15} & \textbf{3.71E-15} & \textbf{3.71E-15} & \textbf{3.71E-15} & 4.00E-15 & 4.00E-15 & 3.85E-15 & \textbf{3.71E-15} \\
    	& $D=15$ & 4.14E-15 & \textbf{4.00E-15} & \textbf{4.00E-15} & \textbf{4.00E-15} & \textbf{4.00E-15} & \textbf{4.00E-15} & \textbf{4.00E-15} & \textbf{4.00E-15} & \textbf{4.00E-15} & \textbf{4.00E-15} & \textbf{4.00E-15} & \textbf{4.00E-15} \\
    	& $D=20$ &  6.70E-15 & 6.55E-15 & 6.27E-15 & 5.84E-15 & \textbf{5.84E-15} & 5.56E-15 & 5.56E-15 & 5.13E-15 & 7.55E-15 & \textbf{4.00E-15} & \textbf{4.00E-15} & \textbf{4.00E-15} \\
    	& $D=25$ & 7.83E-15 & 7.55E-15 & \textbf{7.41E-15} & \textbf{7.41E-15} & \textbf{7.41E-15} & \textbf{7.41E-15} & \textbf{7.41E-15} & \textbf{7.41E-15} & 7.55E-15 & 7.55E-15 & \textbf{7.41E-15} & \textbf{7.41E-15} \\
    	& $D=30$ & 9.54E-15 & 8.83E-15 & 8.69E-15 & 8.12E-15 & 8.40E-15 & 8.26E-15 & 8.12E-15 & 7.69E-15& \textbf{7.55E-15} & \textbf{7.55E-15} & \textbf{7.55E-15} & \textbf{7.55E-15} \\
    	\hline
    	\multirow{6}{*}{$f_4$} & $D=05$ & 3.39E-03 & 2.42E-03 & 2.35E-03 & 1.62E-03 & 1.66E-03 & 1.54E-03 & 1.30E-03 & 1.24E-03 & \textbf{2.27E-04} & \textbf{3.70E-05} & \textbf{2.01E-05} & \textbf{1.64E-05} \\
    	& $D=10$ & 8.95E-03 & 7.72E-03 & 4.50E-03 & 4.45E-03 & \textbf{5.62E-03}  & 4.91E-03 & 4.21E-03 & 3.60E-03 & 7.40E-03 & \textbf{3.67E-04} & \textbf{2.13E-04} & \textbf{1.89E-04} \\
    	& $D=15$ & 1.18E-03 & 9.37E-04 & 3.96E-04 & 3.00E-04 & 2.98E-04 & 3.70E-06 & 4.96E-07 & 4.07E-07 & \textbf{0.00E+00} & \textbf{0.00E+00} & \textbf{0.00E+00} & \textbf{0.00E+00} \\
    	& $D=20$ & \textbf{0.00E+00} &\textbf{0.00E+00} &\textbf{0.00E+00} &\textbf{0.00E+00} &\textbf{0.00E+00} &\textbf{0.00E+00} &\textbf{0.00E+00} &\textbf{0.00E+00} &\textbf{0.00E+00} &\textbf{0.00E+00} &\textbf{0.00E+00} &\textbf{0.00E+00}\\
    	& $D=25$ & \textbf{0.00E+00} &\textbf{0.00E+00} &\textbf{0.00E+00} &\textbf{0.00E+00} &\textbf{0.00E+00} &\textbf{0.00E+00} &\textbf{0.00E+00} &\textbf{0.00E+00} &\textbf{0.00E+00} &\textbf{0.00E+00} &\textbf{0.00E+00} &\textbf{0.00E+00}\\
    	& $D=30$ & \textbf{0.00E+00} &\textbf{0.00E+00} &\textbf{0.00E+00} &\textbf{0.00E+00} &\textbf{0.00E+00} &\textbf{0.00E+00} &\textbf{0.00E+00} &\textbf{0.00E+00} &\textbf{0.00E+00} &\textbf{0.00E+00} &\textbf{0.00E+00} &\textbf{0.00E+00}\\
    	\hline
    	\multirow{6}{*}{$f_5$} & $D=05$ & \textbf{0.00E+00} &\textbf{0.00E+00} &\textbf{0.00E+00} &\textbf{0.00E+00} &\textbf{0.00E+00} &\textbf{0.00E+00} &\textbf{0.00E+00} &\textbf{0.00E+00} &\textbf{0.00E+00} &\textbf{0.00E+00} &\textbf{0.00E+00} &\textbf{0.00E+00}\\
    	& $D=10$ & \textbf{0.00E+00} &\textbf{0.00E+00} &\textbf{0.00E+00} &\textbf{0.00E+00} &\textbf{0.00E+00} &\textbf{0.00E+00} &\textbf{0.00E+00} &\textbf{0.00E+00} &\textbf{0.00E+00} &\textbf{0.00E+00} &\textbf{0.00E+00} &\textbf{0.00E+00}\\
    	& $D=15$ & \textbf{0.00E+00} &\textbf{0.00E+00} &\textbf{0.00E+00} &\textbf{0.00E+00} &\textbf{0.00E+00} &\textbf{0.00E+00} &\textbf{0.00E+00} &\textbf{0.00E+00} &\textbf{0.00E+00} &\textbf{0.00E+00} &\textbf{0.00E+00} &\textbf{0.00E+00}\\
    	& $D=20$ & 2.08E-02 & 4.68E-03 & 3.31E-03 & 5.74E-05 & \textbf{0.00E+00} & \textbf{0.00E+00} & \textbf{0.00E+00} & \textbf{0.00E+00} & \textbf{0.00E+00} & \textbf{0.00E+00} & \textbf{0.00E+00} & \textbf{0.00E+00} \\
    	& $D=25$ & 5.93E-02 & 3.80E-02 & 2.93E-02 & 8.93E-03 & 4.65E-03 & 1.83E-05 & \textbf{0.00E+00} & \textbf{0.00E+00} & \textbf{0.00E+00} & \textbf{0.00E+00} & \textbf{0.00E+00} & \textbf{0.00E+00} \\
    	& $D=30$ & 7.58E-02 & 4.18E-02 & 3.35E-02 & 9.84E-03 & 1.75E-02 &  1.05E-02 & 4.93E-03 & 8.04E-04 &  \textbf{0.00E+00} & \textbf{0.00E+00} & \textbf{0.00E+00} & \textbf{0.00E+00} \\
    	\hline
    	\multirow{6}{*}{$f_6$} & $D=05$ & \textbf{0.00E+00} &\textbf{0.00E+0}0 &\textbf{0.00E+00} &\textbf{0.00E+00} &\textbf{0.00E+00} &\textbf{0.00E+00} &\textbf{0.00E+00} &\textbf{0.00E+00} &\textbf{0.00E+00} &\textbf{0.00E+00} &\textbf{0.00E+00} &\textbf{0.00E+00} \\
    	& $D=10$ & \textbf{0.00E+00} &\textbf{0.00E+0}0 &\textbf{0.00E+00} &\textbf{0.00E+00} &\textbf{0.00E+00} &\textbf{0.00E+00} &\textbf{0.00E+00} &\textbf{0.00E+00} &\textbf{0.00E+00} &\textbf{0.00E+00} &\textbf{0.00E+00} &\textbf{0.00E+00} \\
    	& $D=15$ & \textbf{0.00E+00} &\textbf{0.00E+0}0 &\textbf{0.00E+00} &\textbf{0.00E+00} &\textbf{0.00E+00} &\textbf{0.00E+00} &\textbf{0.00E+00} &\textbf{0.00E+00} &\textbf{0.00E+00} &\textbf{0.00E+00} &\textbf{0.00E+00} &\textbf{0.00E+00} \\
    	& $D=20$ & \textbf{0.00E+00} &\textbf{0.00E+0}0 &\textbf{0.00E+00} &\textbf{0.00E+00} &\textbf{0.00E+00} &\textbf{0.00E+00} &\textbf{0.00E+00} &\textbf{0.00E+00} &\textbf{0.00E+00} &\textbf{0.00E+00} &\textbf{0.00E+00} &\textbf{0.00E+00} \\
    	& $D=25$ & \textbf{0.00E+00} &\textbf{0.00E+0}0 &\textbf{0.00E+00} &\textbf{0.00E+00} &\textbf{0.00E+00} &\textbf{0.00E+00} &\textbf{0.00E+00} &\textbf{0.00E+00} &\textbf{0.00E+00} &\textbf{0.00E+00} &\textbf{0.00E+00} &\textbf{0.00E+00} \\
    	& $D=30$ & \textbf{0.00E+00} &\textbf{0.00E+0}0 &\textbf{0.00E+00} &\textbf{0.00E+00} &\textbf{0.00E+00} &\textbf{0.00E+00} &\textbf{0.00E+00} &\textbf{0.00E+00} &\textbf{0.00E+00} &\textbf{0.00E+00} &\textbf{0.00E+00} &\textbf{0.00E+00} \\
    	\hline
    	\multirow{6}{*}{$f_7$} & $D=05$ & \textbf{0.00E+00} &\textbf{0.00E+0}0 &\textbf{0.00E+00} &\textbf{0.00E+00} &\textbf{0.00E+00} &\textbf{0.00E+00} &\textbf{0.00E+00} &\textbf{0.00E+00} &\textbf{0.00E+00} &\textbf{0.00E+00} &\textbf{0.00E+00} &\textbf{0.00E+00} \\
    	& $D=10$ & \textbf{0.00E+00} &\textbf{0.00E+0}0 &\textbf{0.00E+00} &\textbf{0.00E+00} &\textbf{0.00E+00} &\textbf{0.00E+00} &\textbf{0.00E+00} &\textbf{0.00E+00} &\textbf{0.00E+00} &\textbf{0.00E+00} &\textbf{0.00E+00} &\textbf{0.00E+00} \\
    	& $D=15$ & \textbf{0.00E+00} &\textbf{0.00E+0}0 &\textbf{0.00E+00} &\textbf{0.00E+00} &\textbf{0.00E+00} &\textbf{0.00E+00} &\textbf{0.00E+00} &\textbf{0.00E+00} &\textbf{0.00E+00} &\textbf{0.00E+00} &\textbf{0.00E+00} &\textbf{0.00E+00} \\
    	& $D=20$ & \textbf{0.00E+00} &\textbf{0.00E+0}0 &\textbf{0.00E+00} &\textbf{0.00E+00} &\textbf{0.00E+00} &\textbf{0.00E+00} &\textbf{0.00E+00} &\textbf{0.00E+00} &\textbf{0.00E+00} &\textbf{0.00E+00} &\textbf{0.00E+00} &\textbf{0.00E+00} \\
    	& $D=25$ & \textbf{0.00E+00} &\textbf{0.00E+0}0 &\textbf{0.00E+00} &\textbf{0.00E+00} &\textbf{0.00E+00} &\textbf{0.00E+00} &\textbf{0.00E+00} &\textbf{0.00E+00} &\textbf{0.00E+00} &\textbf{0.00E+00} &\textbf{0.00E+00} &\textbf{0.00E+00} \\
    	& $D=30$ & \textbf{0.00E+00} &\textbf{0.00E+0}0 &\textbf{0.00E+00} &\textbf{0.00E+00} &\textbf{0.00E+00} &\textbf{0.00E+00} &\textbf{0.00E+00} &\textbf{0.00E+00} &\textbf{0.00E+00} &\textbf{0.00E+00} &\textbf{0.00E+00} &\textbf{0.00E+00} \\
    \hline
    \end{tabular}

\end{sidewaystable}

\begin{sidewaystable}
\caption{Statistical mean solution of the R-CSA, B-CSA and PO-CSA for functions $f_8-f_{14}$ with $D \in \{5, 10, 15, 20, 25, 30\}$ and the number of optimizer equals $D$. The stopping criterion goes from  $1\times10^6$ to $8\times10^6$ function evaluations per optimizer. The results are averaged across 25 runs. The best results for each dimension in a given number of function evaluations are highlighted in bold.}\label{tab:quality_solution2}
\tiny
\begin{tabular}{c|c|cccc|cccc|cccc}
	\hline
		\multicolumn{2}{c|}{Algorithm}	& \multicolumn{4}{c|}{R-CSA} & \multicolumn{4}{c|}{B-CSA} & \multicolumn{4}{c}{PO-CSA} \\
		Function & Dimension & $1\times10^6$ & $2\times10^6$ & $4\times10^6$ & $8\times10^6$	& $1\times10^6$ & $2\times10^6$ & $4\times10^6$ & $8\times10^6$	& $1\times10^6$ & $2\times10^6$ & $4\times10^6$ & $8\times10^6$	\\
  	    \hline
    	\multirow{6}{*}{$f_8$} & $D=05$ & \textbf{8.56E-02} & \textbf{8.56E-02} & \textbf{8.56E-02} & \textbf{8.56E-02} & \textbf{8.56E-02} & \textbf{8.56E-02} & \textbf{8.56E-02} & \textbf{8.56E-02} & \textbf{8.56E-02} & \textbf{8.56E-02} & \textbf{8.56E-02} & \textbf{8.56E-02} \\
    	& $D=10$ & 4.12E+01 & 3.98E+01 & \textbf{1.71E-01} & \textbf{1.71E-01} &  \textbf{1.71E-01} & \textbf{1.71E-01} & \textbf{1.71E-01} & \textbf{1.71E-01} & \textbf{1.71E-01} & \textbf{1.71E-01} & \textbf{1.71E-01} & \textbf{1.71E-01} \\
    	& $D=15$ & 9.18E+01 & 9.14E+01 & 8.81E+01 & \textbf{2.57E-01} & \textbf{2.57E-01} & \textbf{2.57E-01} & \textbf{2.57E-01} & \textbf{2.57E-01} & \textbf{2.57E-01} & \textbf{2.57E-01} & \textbf{2.57E-01} & \textbf{2.57E-01} \\
    	& $D=20$ & \textbf{3.43E-01} & 3.43E-01 & \textbf{3.42E-01} & \textbf{3.42E-01} & \textbf{3.43E-01} & \textbf{3.42E-01} & \textbf{3.42E-01} & \textbf{3.42E-01} & \textbf{3.43E-01} & \textbf{3.42E-01} & \textbf{3.42E-01} & \textbf{3.42E-01} \\
    	& $D=25$ & 3.41E+02 & 2.15E+02 & 2.06E+02 & \textbf{4.28E-01} & 9.91E+00 & \textbf{4.28E-01} & \textbf{4.28E-01} & \textbf{4.28E-01} & \textbf{4.28E-01} & \textbf{4.28E-01} & \textbf{4.28E-01} & \textbf{4.28E-01} \\
    	& $D=30$ & 8.51E+02 & 5.57E+02 & 5.49E+02 & 5.03E+01 & 2.69E+02 & \textbf{5.14E-01} & \textbf{5.14E-01} & \textbf{5.14E-01} & \textbf{5.14E-01} & \textbf{5.14E-01} & \textbf{5.14E-01} & \textbf{5.14E-01} \\
    	\hline
    	\multirow{6}{*}{$f_9$} & $D=05$ & 8.70E-16 & 7.28E-16 & 7.28E-16 & 5.86E-16 &  5.86E-16 & 5.86E-16 & 5.86E-16 & \textbf{4.44E-16} &  \textbf{4.44E-16} & \textbf{4.44E-16} & \textbf{4.44E-16} & \textbf{4.44E-16} \\
    	& $D=10$ & 4.14E-15 & 4.00E-15 & 4.00E-15 & 4.00E-15 & \textbf{4.00E-15} & \textbf{3.85E-15} & \textbf{3.85E-15} & \textbf{3.85E-15} & \textbf{4.00E-15} & \textbf{3.85E-15} & \textbf{3.85E-15} & \textbf{3.85E-15} \\
    	& $D=15$ & 7.58E-02 & 5.49E-02 & 4.04E-02 & 8.33E-04 & 3.63E-03 & 4.67E-04 & 3.60E-04 & 2.48E-04 & \textbf{1.64E-13} & \textbf{1.21E-13} & \textbf{1.21E-13} & \textbf{9.99E-14} \\
    	& $D=20$ & 5.53E-01 & 5.30E-01 & 4.58E-01 & 3.25E-01 & 3.58E-01 & 2.54E-01 & 2.26E-01 & 1.24E-01 & \textbf{7.23E-02} & \textbf{2.53E-02} & \textbf{1.75E-02} & \textbf{1.31E-02} \\
    	& $D=25$ & 1.19E+00 & 9.33E-01 & 9.14E-01 & 9.02E-01 & \textbf{1.04E+00} & \textbf{8.20E-01} & 7.53E-01 & 6.64E-01 & 1.16E+00 & 9.28E-01 & \textbf{7.40E-01} & \textbf{6.53E-01} \\
    	& $D=30$ & 2.21E+00 & 2.05E+00 & 1.57E+00 &  1.53E+00 & \textbf{1.47E+00} & 1.20E+00 & 1.11E+00 & 1.03E+00 & 1.78E+00 & \textbf{1.17E+00} & \textbf{1.10E+00} & \textbf{9.41E-01} \\
    	\hline
    	\multirow{6}{*}{$f_{10}$} & $D=05$ & 1.88E-02 & 1.80E-02 & 1.76E-02 & 1.71E-02 & \textbf{1.71E-02} & 1.60E-02 & \textbf{1.45E-02} & 1.36E-02 & 1.73E-02 & \textbf{1.49E-02} & \textbf{1.45E-02} & \textbf{1.25E-02} \\
    	& $D=10$ & 3.58E-02 & 2.93E-02 &  2.91E-02 & 2.65E-02 & \textbf{2.88E-02} & \textbf{2.65E-02} & 2.51E-02 & 1.95E-02 & 3.40E-02 & \textbf{2.65E-02} & \textbf{2.50E-02} & \textbf{1.81E-02} \\
    	& $D=15$ & 1.88E-02 & 1.71E-02 & 1.60E-02 & 1.46E-02 & 1.45E-02 & 1.40E-02 & 1.23E-02 & \textbf{1.00E-02} & \textbf{1.40E-02} & \textbf{1.31E-02} & \textbf{1.22E-02} & \textbf{1.00E-02} \\
    	& $D=20$ & 1.56E-02 & 1.26E-02 & 9.21E-03 & 8.80E-03 & \textbf{1.17E-02} & 8.61E-03 & 8.60E-03 & 8.33E-03 & \textbf{1.17E-02 }& \textbf{8.59E-03} & \textbf{7.97E-03} & \textbf{7.41E-03} \\
    	& $D=25$ & 1.26E-02 & 9.34E-03 & 7.93E-03 & 5.43E-03 & \textbf{8.17E-03} & \textbf{5.92E-03} & \textbf{5.43E-03} & 5.34E-03 &  9.51E-03 & 8.85E-03 & 6.48E-03 & \textbf{1.54E-03}  \\
    	& $D=30$ & 5.63E-03 & 4.96E-03 & 3.18E-03 & 2.89E-03 & 3.78E-03 & 3.28E-03 & 3.18E-03 & 2.89E-03 & \textbf{3.41E-03} & \textbf{2.05E-03} & \textbf{1.85E-03} & \textbf{1.14E-03} \\
    	\hline
    	\multirow{6}{*}{$f_{11}$} & $D=05$ & 3.47E-02 & 2.82E-02 & 1.14E-02 & 2.57E-03 & 1.19E-02 & 1.11E-02 & 6.54E-03 & 2.41E-03 & \textbf{4.77E-05} & \textbf{1.33E-05} & \textbf{3.57E-07} & \textbf{1.93E-07}\\
    	& $D=10$ & 5.61E-01 & 4.81E-01 & 4.60E-01 & 2.47E-01 & 3.67E-01 & 3.23E-01 & 3.08E-01 & 2.37E-01 & \textbf{1.86E-01} & \textbf{1.52E-01} & \textbf{7.77E-02} & \textbf{7.40E-02}\\
    	& $D=15$ & 1.33E+00 & 1.24E+00 & 1.16E+00 & 1.10E+00 & 1.17E+00 & 1.05E+00 & 9.84E-01 & 9.25E-01 & \textbf{1.15E+00} & \textbf{1.00E+00} & \textbf{9.79E-01} & \textbf{9.15E-01} \\
    	& $D=20$ & 2.95E+00 & 2.93E+00 & 2.85E+00 & 2.62E+00 & 2.53E+00 & 2.50E+00 & 2.44E+00 & 2.17E+00 & \textbf{2.50E+00} & \textbf{2.41E+00} & \textbf{2.23E+00} & \textbf{2.11E+00} \\
    	& $D=25$ & 4.94E+00 & 4.60E+00 & 4.48E+00 & 4.41E+00 & 4.42E+00 & 4.16E+00 & 4.13E+00 & 3.81E+00 & \textbf{4.41E+00} & \textbf{4.12E+00} & \textbf{4.11E+00} & \textbf{3.67E+00} \\
    	& $D=30$ & 7.34E+00 & 7.26E+00 & 6.91E+00 & 6.34E+00 & \textbf{6.39E+00} & \textbf{6.34E+00} & \textbf{6.21E+00} & 6.14E+00 & 7.27E+00 & 7.18E+00 & 6.74E+00 & \textbf{6.11E+00 }\\
    	\hline
    	\multirow{6}{*}{$f_{12}$} & $D=05$ & 1.19E-01 & 7.97E-02 & 4.08E-02 & 4.01E-02 & 6.73E-04 & 9.51E-05 & 8.88E-05 & 7.48E-05 & \textbf{6.21E-08} & \textbf{0.00E+00} & \textbf{0.00E+00} & \textbf{0.00E+00}\\
    	& $D=10$ & 3.28E+00 & 3.03E+00 & 2.95E+00 & 2.76E+00 & 2.87E+00 & 2.83E+00 & 2.71E+00 & 2.63E+00 & \textbf{2.85E+00} & \textbf{2.80E+00} & \textbf{2.70E+00} & \textbf{2.61E+00} \\
    	& $D=15$ & 9.20E+00 & 9.12E+00 & 8.64E+00 & 7.36E+00 & 8.15E+00 & 8.07E+00 & 7.78E+00 & 7.23E+00 & \textbf{8.05E+00} & \textbf{8.03E+00} & \textbf{7.75E+00} & \textbf{7.22E+00 }\\
    	& $D=20$ & 1.94E+01 & 1.73E+01 & 1.68E+01 & 1.55E+01 & 1.52E+01 & 1.52E+01 & 1.50E+01 & 1.46E+01 & \textbf{1.50E+01} & \textbf{1.49E+01} & \textbf{1.48E+01} & \textbf{1.43E+01 }\\
    	& $D=25$ & 2.18E+01 & 2.15E+01 & 2.12E+01 & 2.04E+01 & 2.03E+01 & 2.02E+01 & \textbf{1.92E+01} & \textbf{1.83E+0}1 & \textbf{2.02E+01} & \textbf{2.00E+01} & \textbf{1.92E+01} & \textbf{1.83E+01 }\\
    	& $D=30$ & 3.20E+01 & 3.19E+01 & 3.02E+01 & 3.01E+01 & \textbf{2.94E+01} & \textbf{2.94E+01} & \textbf{2.91E+01} & \textbf{2.88E+01} & 3.19E+01 & 3.02E+01 & 3.00E+01 & 2.99E+01 \\
    	\hline
    	\multirow{6}{*}{$f_{13}$} & $D=05$ & 2.25E-01 & 8.01E-02 & 6.21E-02 & 2.21E-03 & 4.00E-02 & 4.00E-02 & 2.90E-03 & 1.29E-03 & \textbf{1.32E-06} & \textbf{4.61E-08} & \textbf{1.36E-08} & \textbf{1.23E-08}\\
    	& $D=10$ & 2.53E+00 & 2.36E+00 & 2.32E+00 & 2.00E+00 & \textbf{2.05E+00} & 2.05E+00 & 1.96E+00 & \textbf{1.85E+00} & 2.47E+00 & \textbf{2.01E+00} & \textbf{1.88E+00} & \textbf{1.85E+00} \\
    	& $D=15$ & 5.61E+00 & 6.57E+00 & 5.49E+00 & 4.89E+00 & 5.07E+00 & 5.06E+00 & 5.00E+00 & \textbf{4.84E+00} & \textbf{5.00E+00} & \textbf{4.94E+00} & \textbf{4.85E+00} & \textbf{4.84E+00 }\\
    	& $D=20$ & 1.35E+01 & 1.15E+01 & 1.06E+01 & 1.05E+01 & \textbf{1.02E+01} & \textbf{1.02E+01} & \textbf{1.01E+01} & 9.62E+00 & 1.20E+01 & 1.10E+01 & \textbf{1.01E+01} & \textbf{9.60E+00} \\
    	& $D=25$ & 1.73E+01 & 1.38E+01 & 1.36E+01 & 1.26E+01 & \textbf{1.41E+01} & \textbf{1.38E+01} & \textbf{1.36E+01} & 1.25E+01 & 1.50E+01 & \textbf{1.38E+01} & \textbf{1.36E+01} & \textbf{1.24E+01} \\
    	& $D=30$ & 2.31E+01 & 2.26E+01 & 2.20E+01 & 2.08E+01 & 2.12E+01 & 2.06E+01 & 2.04E+01 & \textbf{1.92E+01} & \textbf{2.10E+01} & \textbf{2.00E+01} & \textbf{2.00E+01} & 2.00E+01 \\
    	\hline
		\multirow{6}{*}{$f_{14}$} & $D=05$ & 9.02E-02 & 8.66E-02 & 8.63E-02 & 8.57E-02 & 8.60E-02 & 8.60E-02 & 8.59E-02 & \textbf{8.56E-02} & \textbf{8.57E-02} & \textbf{8.57E-02} & \textbf{8.56E-02} & \textbf{8.56E-02} \\
    	& $D=10$ & 9.62E+01 & 9.12E+01 & 5.33E+01 & 4.17E+01 & 5.25E+01 & 4.32E+01 & 4.00E+01 & 3.35E+01 & \textbf{3.27E+00} & \textbf{2.26E+00} & \textbf{1.44E+00} & \textbf{1.31E+00} \\
    	& $D=15$ & 7.33E+02 & 6.94E+02 & 6.93E+02 & 6.78E+02 & \textbf{6.48E+02} & 6.46E+02 & 6.29E+02 & 6.12E+02 & \textbf{6.48E+02} & \textbf{6.44E+02} & \textbf{6.19E+02} & \textbf{5.85E+02} \\
    	& $D=20$ & 1.69E+03 & 1.67E+03 & 1.56E+03 & 1.53E+03 & \textbf{1.54E+03} & \textbf{1.52E+03} & \textbf{1.50E+03} & \textbf{1.50E+03} & 1.61E+03 & 1.60E+03 & \textbf{1.50E+03} & \textbf{1.50E+03} \\
    	& $D=25$ & 2.17E+03 & 2.15E+03 & 1.97E+03 & 1.97E+03 & \textbf{2.00E+03} & \textbf{1.98E+03} & \textbf{1.96E+0}3 & 1.95E+03 & 2.08E+03 & 2.01E+03 & \textbf{1.96E+03} & \textbf{1.93E+03} \\
    	& $D=30$ & 3.18E+03 & 3.15E+03 & 3.03E+03 & 3.00E+03 & \textbf{3.00E+03} & \textbf{3.00E+03} & \textbf{2.97E+03} & \textbf{2.96E+03} & 3.14E+03 & \textbf{3.00E+03} & \textbf{2.97E+03} & \textbf{2.96E+03} \\

    \hline
    \end{tabular}

\end{sidewaystable}

The results of the B-CSA are always better than or equal to the ones from the R-CSA, as expected, given that the initial values of B-CSA were chosen after exhaustive executions with several different initial generation temperatures. Since the PO-CSA is parameter-free, the results of the R-CSA and B-CSA, then, represents the minimum goal and the golden goal for the PO-CSA, respectively. No so surprisingly, the PO-CSA reached equal or better solutions than the R-CSA in $98.21\%$ of the cases, only losing in $f_2$ in 5 cases and 1 case in $f_3$. However, the effectiveness of the PO-CSA exceeds expectations when it is equal or better than the B-CSA in $86.90\%$ of the tests. The remaining $13.10\%$ cases, where PO-CSA looses for B-CSA, are distributed in functions $f_2$, $f_3$, $f_4$, $f_9$, $f_{10}$, $f_{11}$, $f_{12}$, $f_{13}$, and $f_{14}$. For the PO-CSA, the worst case is for the function $f_2$. 

\subsubsection{The Resilience of the PO-CSA}
\label{subsubsec:analysis_f2}

In this set of experiments, we compared the quality of the solutions found by the R-CSA, B-CSA and PO-CSA for $f_2$, whose PO-CSA's failure rate was the highest in the previous set of experiment, to understand the behavior of the three algorithms when we have an increased budget of function evaluations. In order to perform these experiments, we used dimension $D \in \{5,10,15,20,25,30\}$ and the number of optimizers equals $D$. The stopping criterion for the three algorithms was varied from $1.0\times10^6$ to $16.0\times10^6$ and results were averaged over 25 runs.

\begin{figure*}[!ht]
	\centering
	\psfragscanon
	\psfrag{d}[][][0.8]{Dimension} 
	\psfrag{v}[][][0.8]{Energy} 
	\psfrag{R-CSA}[][][0.3]{R-CSA}
	\psfrag{B-CSA}[][][0.3]{B-CSA}
	\psfrag{PO-CSA}[][][0.3]{PO-CSA}
	\subfloat[$1\times10^6$ function evaluations per optimizer\label{fig:pot-csa_f2_1mi}]{
		\includegraphics[width=0.42\columnwidth]{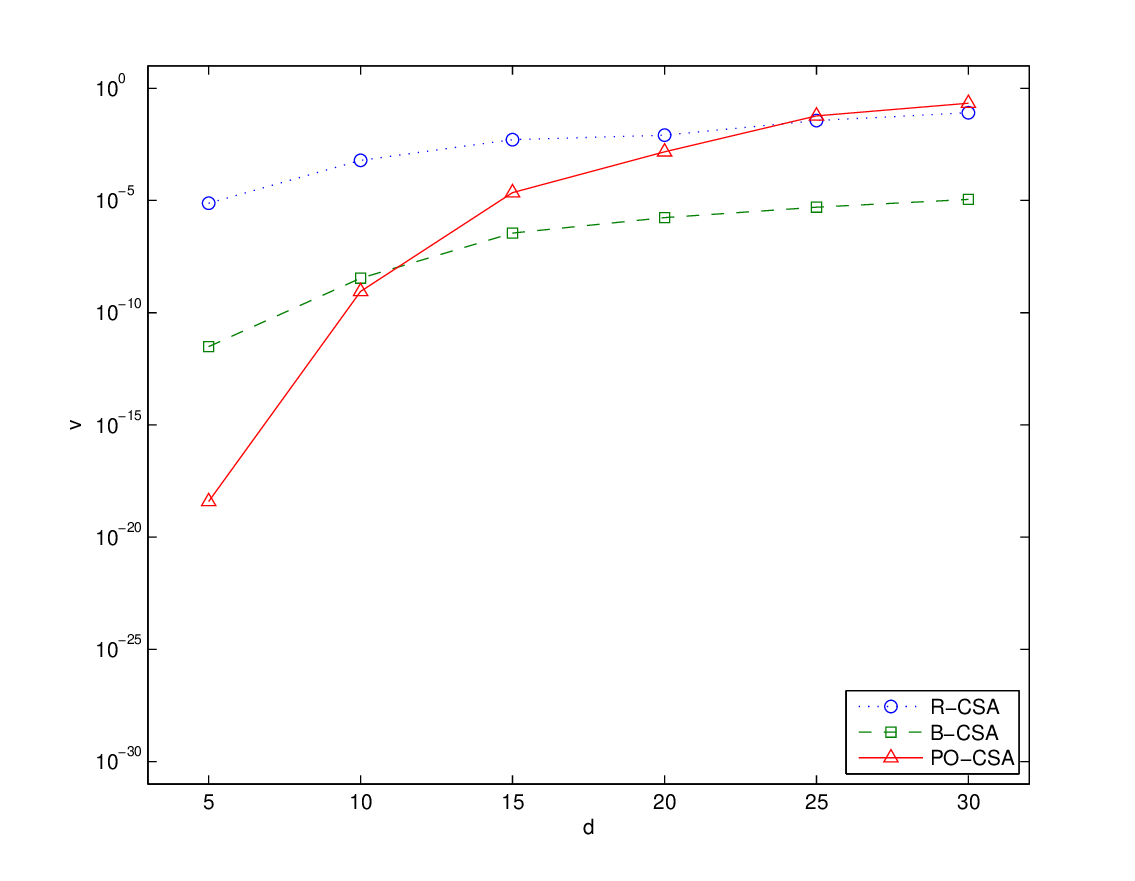}
    }
    \subfloat[$2\times10^6$ function evaluations  per optimizer\label{fig:pot-csa_f2_2mi}]{
		\includegraphics[width=0.42\columnwidth]{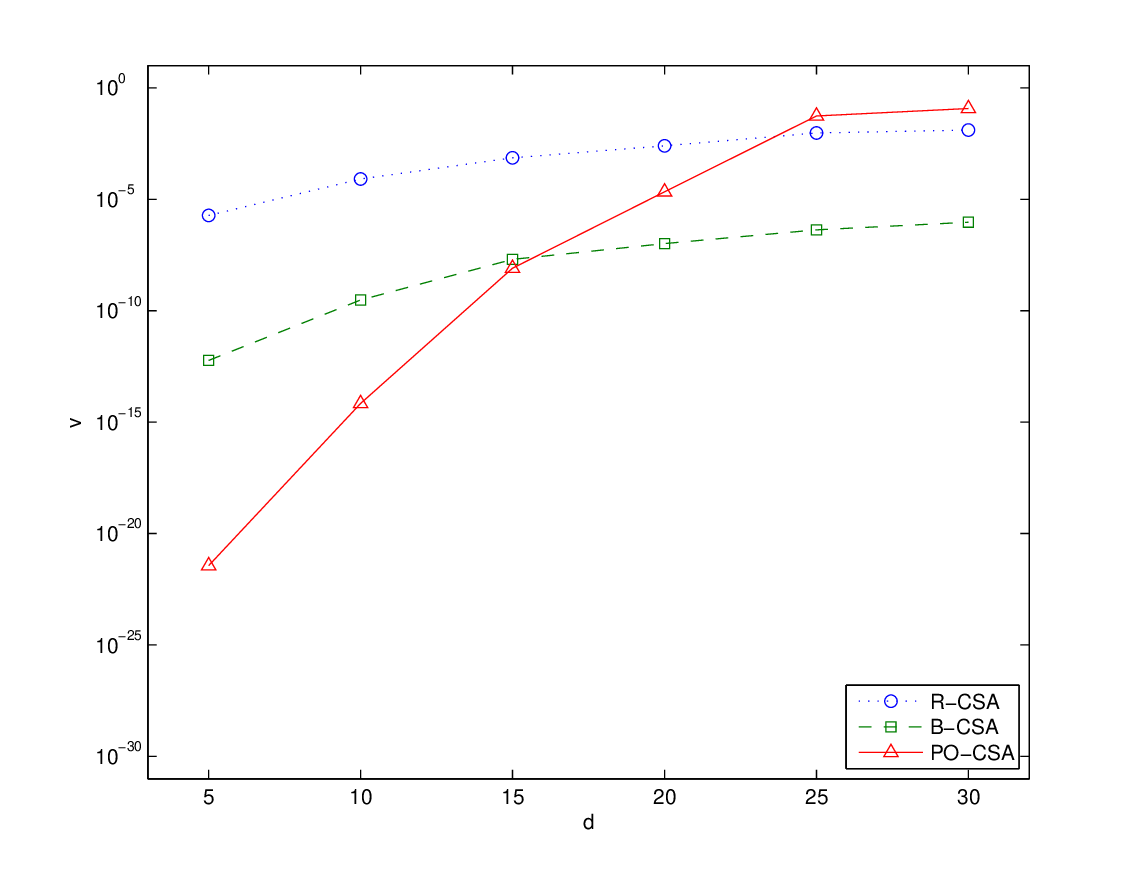}
    }
    
    \subfloat[$4\times10^6$ function evaluations  per optimizer\label{fig:pot-csa_f2_4mi}]{
		\includegraphics[width=0.42\columnwidth]{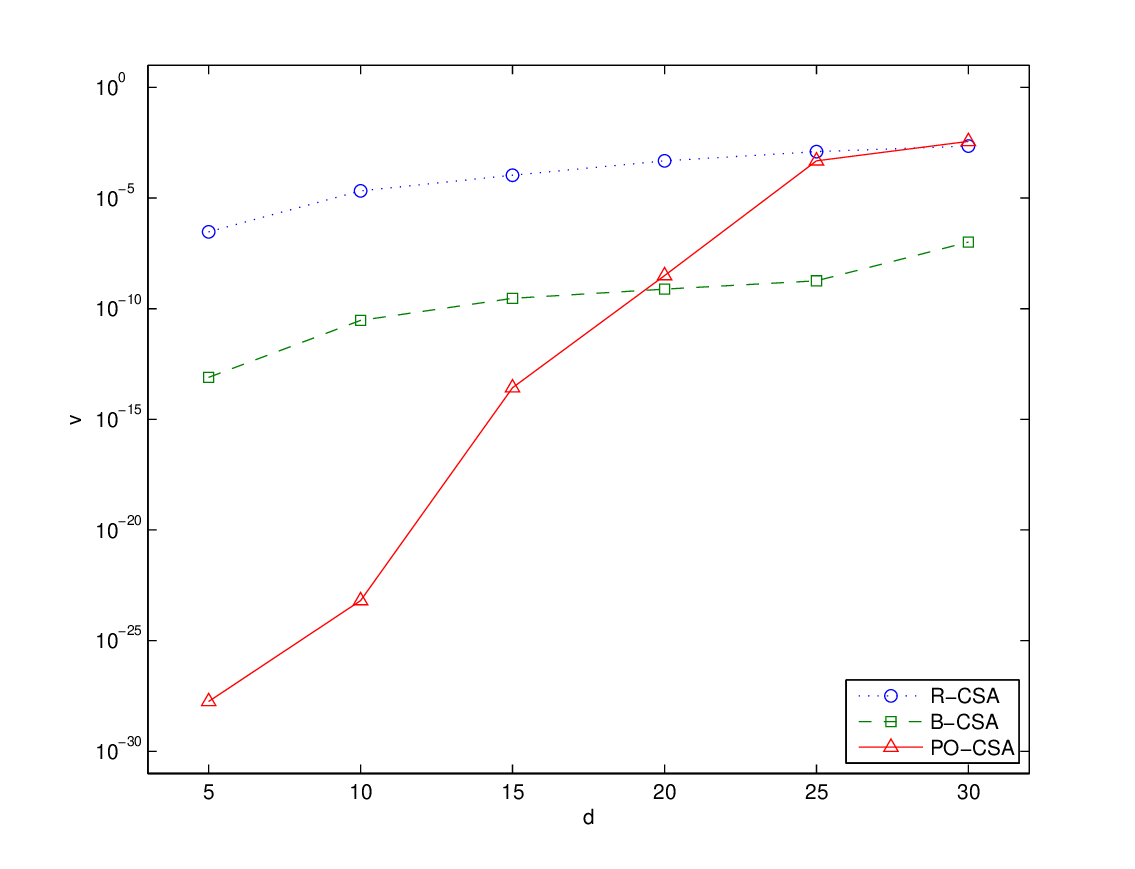}

    }
    \subfloat[$8\times10^6$ function evaluations  per optimizer\label{fig:pot-csa_f2_8mi}]{
		\includegraphics[width=0.42\columnwidth]{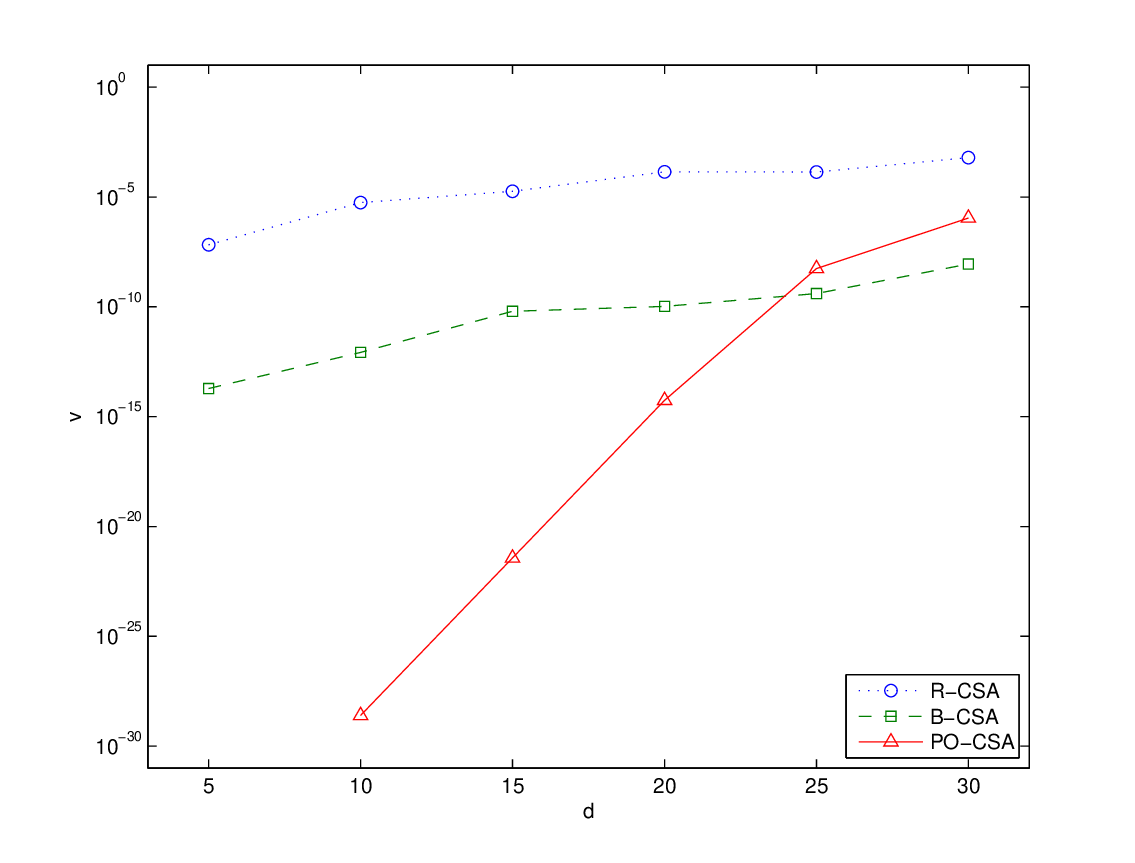}
    }
    
    \subfloat[$12\times10^6$ function evaluations  per optimizer\label{fig:pot-csa_f2_12mi}]{
		\includegraphics[width=0.42\columnwidth]{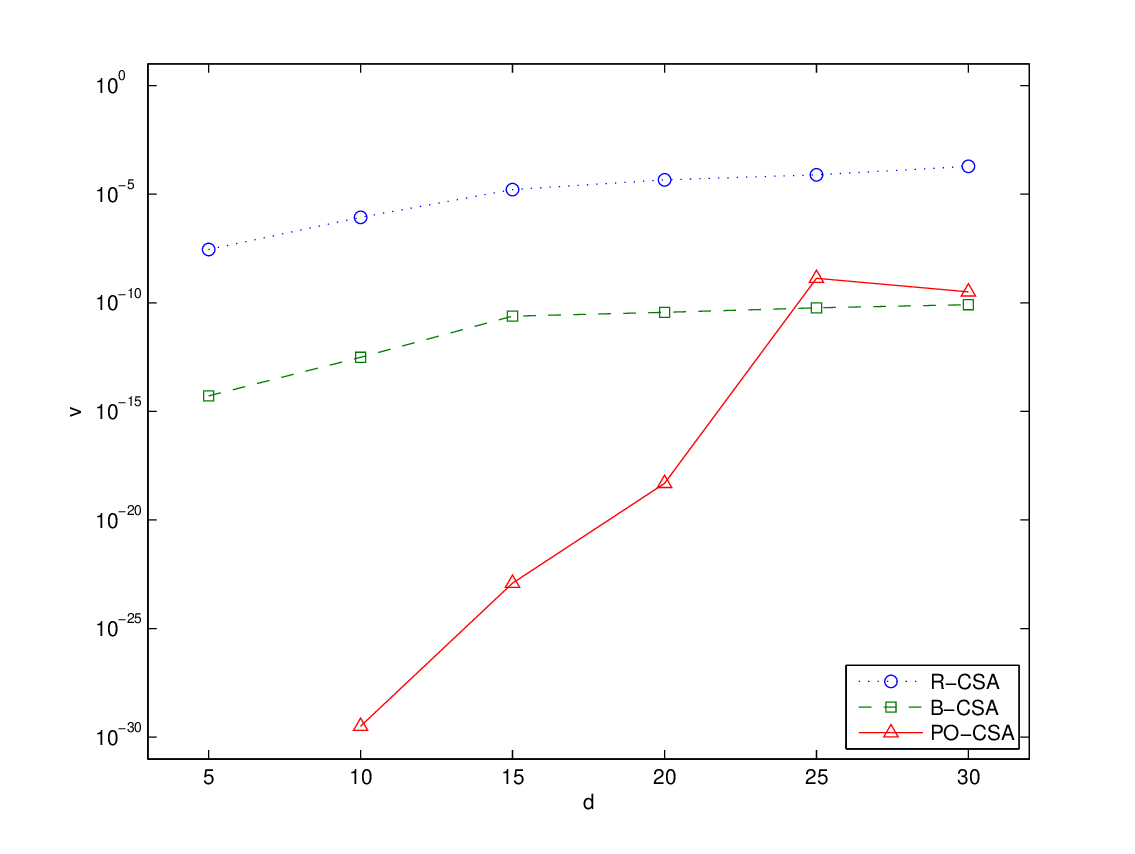}
    }
    \subfloat[$16\times10^6$ function evaluations  per optimizer\label{fig:pot-csa_f2_16mi}]{
		\includegraphics[width=0.42\columnwidth]{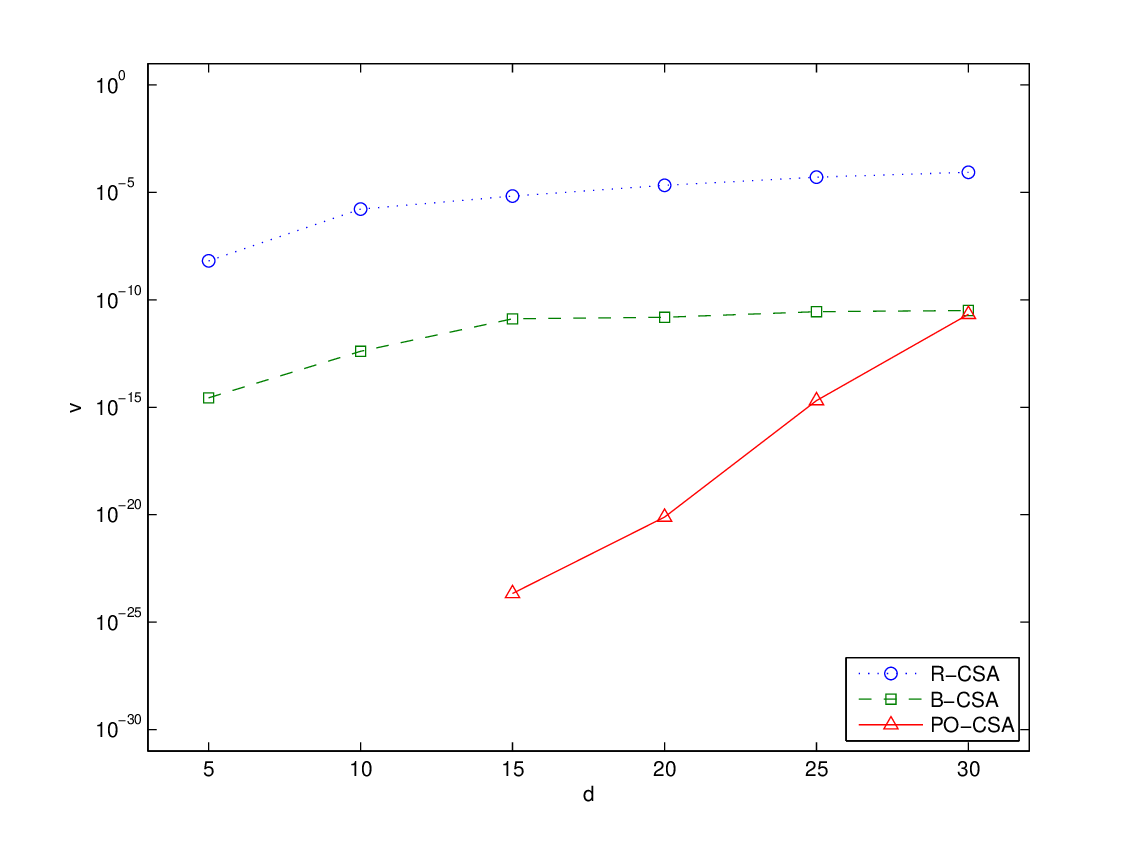}
    }
    \caption{Average energy of $f_2$ over 25 runs for different numbers of function evaluations per optimizer. The number of optimizers is the same as the dimension. In (d) and (e) for $D=5$, and (f) for $D=5,10$, there are no points for PO-CSA because its costs are zero (global optimum) and not defined in logarithmic scale. The larger the number of function evaluations, the better the results of the PO-CSA because it does not stagnates its search.}
  \label{fig:pot-csa_f2}
\end{figure*}

Observe in Fig.~\ref{fig:pot-csa_f2} that, although the PO-CSA does not beat the B-CSA in all cases, even the R-CSA for $D \in \{25,30\}$ with 1 and 2 million of function evaluations, the larger the number of functions evaluations, the better the results of the PO-CSA. For 4, 8 and 12 million of function evaluations, the PO-CSA wins the R-CSA in all cases, but it still does not win the B-CSA for $D \in \{25,30\}$. For 16 million of function evaluations, the PO-CSA wins the B-CSA in 100\% of the cases. Besides, the PO-CSA found the global optimum for $D=5$ when running 8 and 12 million of functions evaluations as well as for $D=5,10$ for 16 million of functions evaluations, which represents the lack of marking dots for PO-CSA in Fig.~\ref{fig:pot-csa_f2_8mi}, Fig.~\ref{fig:pot-csa_f2_12mi} and Fig.~\ref{fig:pot-csa_f2_16mi} because zero is not defined in logarithmic scale. This behavior happens because the PO-CSA does not stagnates its searching processes, since the generation temperature is not monotonically decreasing. On the other hand, increasing the number of function evaluations of the CSA can stagnate its search because the generation temperature can reach a vary small value.

\section{Conclusion}
\label{sec:conclusion}
We have presented an algorithm based on Coupled Simulated Annealing (CSA) that controls all its setting parameters. The Perpetual Orbit Coupled Simulated Annealing (PO-CSA) was designed by applying a novel technique called Perpetual Orbit (PO) to the CSA algorithm. Due to the coupling, the CSA is robust against the initialization of the acceptance temperature. However, the generation temperature requires tuning to avoid search stagnation. The PO technique arises as a solution that control the generation temperature, which no longer decreases monotonically. In principle, PO can also be applied to other ensemble- and population-based algorithms that have a dispersion variable.

We performed four sets of experiments to prove the effectiveness of the PO-CSA. The experiments were performed using 14 objective functions. In the first set of experiments, we compared the quality of solutions found by the PO-CSA to those found by the Cuckoo Search via L\'{e}vi Flight, Differential Evolution, Particle Swarm Optimization and Genetic Algorithm. Results show that the PO-CSA found equal or better solutions in 85.71\% of the cases.

In other two sets of experiments, we compared the PO-CSA to the CSA with two different initialization procedures. The first procedure, called R-CSA, has its initial generation temperature initialized as a random number. The second one, B-CSA, shows the best solutions through an exhaustive search of the initial generation temperature. These two procedures had the same monotonically decreasing temperature schedule. Results from the R-CSA are the first frontiers that the PO-CSA must overcome. Those from the B-CSA express excellence in quality of solution. The comparison shows show that the PO-CSA overcomes the R-CSA in $98.21\%$ of the cases and the B-CSA in $86.90\%$. 

For the worst case, where the PO-CSA performed less than the B-CSA,  the experiments show that the R-CSA and B-CSA do not improve the quality of the solution when the number of function evaluations increases. However, the same does not happen with the PO-CSA due to the non-stagnation of the generation temperature introduced by the PO technique. In this case, sufficiently increasing the number of function evaluations makes the solutions of the PO-CSA equal or better than those of the B-CSA in $100\%$ of the cases. 

Lastly, the set of experiments with the analysis of the generation temperature control introduced by the PO technique  shows why the PO-CSA presents superior performance. By allowing the independent control of its initial setting parameter, the PO technique makes the PO-CSA parameter-free and capable of reaching better solutions than the original CSA in the majority of functions for all competitor methods tested.

\begin{appendices}

\section{The Objective Functions}
\label{appendix}

In order to guarantee the satisfactory of the PO-CSA, experiments were carried out repeatedly in 14 reference functions $f_1-f_{14}$, with $D$ dimensions, used in \cite{samuel2010}. The global optimum for all functions is $f^*(\textbf{x})=0$. They are detailed below.

 \begin{enumerate}
   \item \textbf{(Group 1: $f_1-f_2$)} Unimodal and simple multimodal functions: the function 1 is an easy unimodal sphere function. The second function is the ubiquitous Rosenbrock’s function, which is very often used for testing optimization algorithms. They are described in Table~\ref{tab:function_table1}.

\begin{table}[!ht] 
	\centering
	\caption{Unimodal and simple multimodal functions: test problems, Group 1. The minimum of both functions is zero. The dimensionality of These problems can be adjusted with the term $D$.} 
		\begin{tabular}{ccc} 
		\hline 
		No. 	& 	Function $f(\textbf{x)}$ & Input range\\
		\hline
			1 & $f_1=\sum \limits_{i=1}^D x_i^2$ & [-100,100] \\
			2 & $f_2=\sum \limits_{i=1}^{D-1} (1-x_i)^2 + 100(x_{i+1}-x_i^2)^2$ & [-2.048,2.048]\\
		\hline
	\end{tabular}
	\label{tab:function_table1}
\end{table}

   \item \textbf{(Group 2: $f_3-f_8$)} Multimodal Functions: A collection of multidimensional and multimodal continuous functions was chosen from the literature to be used as test cases. These functions feature many local minima and therefore are regarded as being
difficult to optimize~\cite{tu2004robust,yao1999evolutionary}. They are detailed in Table~\ref{tab:function_table2}.

\begin{table}[!ht] 
\scriptsize
	\centering
	\caption{Multimodal functions: test problems, group 2. These functions present many local minima and are considered hard problems to optimize, particularly in large dimensions. The minimum of all these functions is zero. The dimensionality of these functions can be adjusted with the term $D$.} 
	\begin{tabular}{ccc} 
		\hline 
		No. 	& 	Function $f(\textbf{x)}$ & Input range\\
		\hline
			3 & $f_3 =  -20exp\left[ -0.2\sqrt{\frac{1}{D}\sum \limits_{i=1}^D x_i^2} \right] - exp\left[ \frac{1}{D} \sum \limits_{i=1}^D cos(2 \pi x_i) \right] + 20 + e  $ & [-32.768,32.768] \\
			4 & $f_4 = \sum \limits_{i=1}^D \frac{x_i^2}{4000} - \prod \limits_{i=1}^D cos \left(\frac{x_i}{\sqrt{i}}\right) + 1$ & [-600,600]\\
			5 & $f_5 = \sum \limits_{i=1}^D \left \{  \sum \limits_{k=0}^{20} \left[ (0.5)^k cos \left(2 \pi 3^k (x_i + 0.5)\right) \right] \right \}  - D \sum \limits_{k=0}^{20} \left[0.5^k cos ( \pi 3^k ) \right]  $ & [-0.5,0.5]\\
			6 & $f_6 = \sum \limits_{i=1}^D x_i^2 - 10cos(2 \pi x_i) + 10 $ & [-5.12,5.12] \\
			7 & $f_7 = \sum \limits_{i=1}^D y_i^2 - 10cos(2 \pi y_i) + 10$ , $y_i = 
					\left\{
					\begin{array}{rc}
						x_i \quad |x_i| < \frac{1}{2} \\
						\frac{round(2x_i)}{2} \quad |x_i| \ge \frac{1}{2}
					\end{array}
					,~ \forall i=1~,~\ldots~,~D
					\right. 		
					
					$& [-5.12,5.12]\\
			8 & $f_8 = 419 \times D + \sum \limits_{i=1}^D x_i sin (|x_i|^{\frac{1}{2}}) $ & [-500,500] \\
		\hline
	\end{tabular}
	\label{tab:function_table2}
\end{table}
   
   \item \textbf{(Group 3: $f_9-f_{14}$)} Nonseparable Functions: The functions of the Group 2 are considered hard to optimize, but they can possibly be separable. This condition means that the minimization problem can be solved using $D$ unidimensional searches, where $D$ is the dimension of the problem. Various real-life optimization problems are nonseparable. Therefore, to approximate these problems, in this group of test problems, we use a set of rotated versions of the Group 2. The rotated functions preserve the same shape characteristics as those of the original functions but cannot be solved by $D$ unidimensional searches. As described in~\cite{samuel2010}, in order to rotate a function, we multiply the argument $\textbf{x}$ by an orthogonal rotation matrix $\textbf{M}$ to obtain a new argument $\textbf{z}$ for the rotated function. This rotation matrix was obtained by using Salomon's method~\cite{salomon1996re}. So, we can define the first 5 functions in this group

\begin{equation}
	f_n(\textbf{x}) = f_{n-6}(\textbf{z}),~~~\forall n=9,\ldots,13,
\end{equation}   
with $\textbf{z}=\textbf{Mx}$. The last function $f_{14}$ is defined as follows:

\begin{equation}
	f_{14}(\textbf{x}) = f_{8}(\textbf{z})
\end{equation}   
with 

\begin{align}
z_i = 
		\left\{ 
					\begin{array}{rc}
					y_i sin\left( |y_i|^{\frac{1}{2}} \right), \quad |y_i| \le 500\\
					0.001 (|y_i|-500)^2, \quad |y_i| > 500
					\end{array}
					,
		\right.\\
		\forall i=1,\ldots,D \nonumber
\end{align}
\begin{equation}
\textbf{y} = \textbf{y}' + 420.96 
\end{equation} 
\begin{equation}
\textbf{y}' = \textbf{M}(\textbf{x}-420.96). 
\end{equation}

This condition is necessary to keep the global optimum of the original Schwefel’s function, which is located at $[420.96,420.96,\ldots,420.96]$, within the search range after rotation.

 \end{enumerate}

\section*{Acknowledgment}

This research was supported by the High-Performance Computing Center at UFRN (NPAD/UFRN).




\end{appendices}


\bibliography{sn-bibliography}

\end{document}